\documentclass[a4paper,11pt]{article}
\usepackage{bbm}
\usepackage{amsfonts}
\usepackage{amssymb}
\usepackage{amsthm}
\usepackage{amsmath}
\usepackage{graphicx}
\usepackage{empheq}
\usepackage{float}
\usepackage{indentfirst}
\usepackage{mathrsfs}
\usepackage{multirow}
\usepackage{booktabs}
\usepackage{graphics}
\usepackage[percent]{overpic}
\usepackage{color}
\usepackage{epstopdf}
\usepackage{authblk}
\usepackage[justification=centering]{caption}

\graphicspath{ {./figs_Num_CHSD/figs_Num_CHSD/} }

\allowdisplaybreaks[4]

\newtheorem{mydef}{Definition}
\newtheorem{mytheo}{Theorem}
\newtheorem{myassum}{Assumption}

\newtheorem{lemma}{\textbf{Lemma}}
\newtheorem{corollary}{\textbf{Corollary}}
\newtheorem{remark}{\textbf{Remark}}

\newcommand{\be}{\begin{equation}}
\newcommand{\ee}{\end{equation}}
\newcommand{\bes}{\begin{equation*}}
\newcommand{\ees}{\end{equation*}}
\newcommand{\pp}[2]{\frac{\partial{#1}}{\partial{#2}}}

\newcommand{\ignore}[1]{}

\def\bea{\begin{eqnarray}}
\def\eea{\end{eqnarray}}
\def\bt{\begin{theorem}}
\def\et{\end{theorem}}
\def\bl{\begin{lemma}}
\def\el{\end{lemma}}
\def\br{\begin{remark}}
\def\er{\end{remark}}
\def\bp{\begin{proposition}}
\def\ep{\end{proposition}}
\def\bc{\begin{corollary}}
\def\ec{\end{corollary}}
\def\bd{\begin{definition}}
\def\ed{\end{definition}}

\def\vp{\varphi}

\def\non{\nonumber }
\def\ub{\mathbf{u}}
\def\vb{\mathbf{v}}
\def\baru{\overline{\mathbf{u}}}
\def\btau{\boldsymbol{\tau}}
\def\BH{\mathbf{H}}

\def\ar1{\varphi_h^{k+1}}
\def\var0{\varphi_h^{k}}
\def\bu1{\overline{\mathbf{u}}_h^{k+1}}

\begin{document}

\title{Error estimate of a decoupled numerical scheme for the Cahn-Hilliard-Stokes-Darcy system}

\author[1]{Wenbin Chen\thanks{ wbchen@fudan.edu.cn}}
\author[2]{Daozhi Han\thanks{  handaoz@mst.edu}}
\author[3]{ Cheng Wang\thanks{  cwang1@umassd.edu}}
\author[1]{Shufen Wang\thanks{ 17110180015@fudan.edu.cn}}
\author[4]{Xiaoming Wang\thanks{ wxm.math@outlook.com}}
\author[1]{Yichao Zhang \thanks{yichaozhang16@fudan.edu.cn}}

\affil[1]{School of Mathematical Sciences, Fudan University,  Shanghai, China}
\affil[2]{Department of Mathematics and Statistics, Missouri University of Science and Technology, Rolla, MO}
\affil[3]{Department of Mathematics, University of Massachusetts Dartmouth, North Dartmouth, MA}
\affil[4]{Department of Mathematics, Southern University of Science and Technology, Shenzhen, China }

\maketitle

\begin{abstract}
	{We analyze a fully discrete finite element numerical scheme for the Cahn-Hilliard-Stokes-Darcy system that models two-phase flows in coupled free flow and porous media. To avoid a well-known difficulty associated with the coupling between the Cahn-Hilliard equation and the fluid motion, we make use of the operator-splitting in the numerical scheme, so that these two solvers are decoupled, which in turn would greatly improve the computational efficiency. The unique solvability and the energy stability have been proved in~\cite{CHW2017}. In this work, we carry out a detailed convergence analysis and error estimate for the fully discrete finite element scheme, so that the optimal rate convergence order is established in the energy norm, i.e.,, in the $\ell^\infty (0, T; H^1) \cap \ell^2 (0, T; H^2)$ norm for the phase variables, as well as in the $\ell^\infty (0, T; H^1) \cap \ell^2 (0, T; H^2)$ norm for the velocity variable. Such an energy norm error estimate leads to a cancellation of a nonlinear error term associated with the convection part, which turns out to be a key step to pass through the analysis. In addition, a discrete $\ell^2 (0;T; H^3)$ bound of the numerical solution for the phase variables plays an important role in the error estimate, which is accomplished via a discrete version of Gagliardo-Nirenberg inequality in the finite element setting.} 

\end{abstract}



\section{Introduction}
\begin{figure}[ht]
  \begin{center}
    \includegraphics[width=0.7\textwidth]{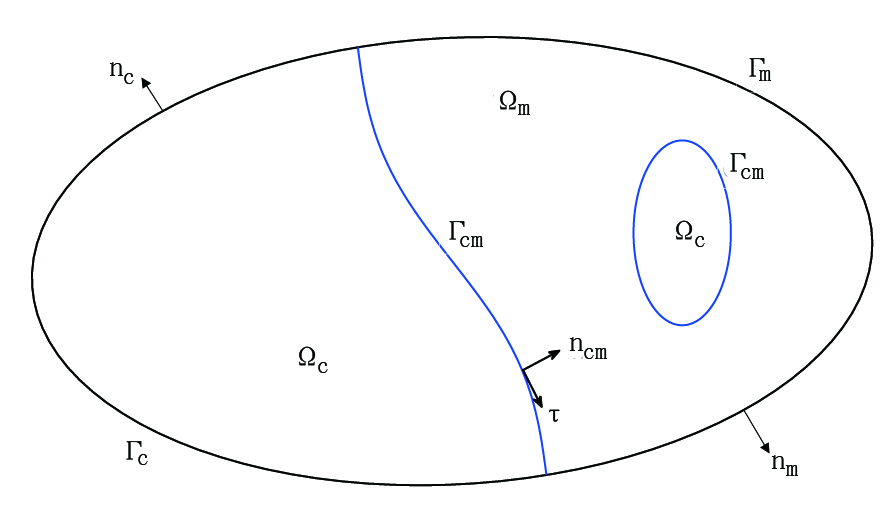}
  \end{center}
  \caption{Schematic illustration of karst geometry in 2D}
  \label{domain}
\end{figure}

In many applications such as contaminant transport in karst aquifer, oil recovery in karst oil reservoir, proton exchange membrane fuel cell technology, cardiovascular modeling, multiphase flows in conduit and in porous media interact with each other, and therefore have to be considered together. Geometric configurations that consist of both conduit and porous media are termed as karstic geometry.
In this article we aim to analyze a decoupled numerical algorithm for solving the Cahn-Hilliard-Stokes-Darcy model (CHSD) for two-phase flows in Karst geometry--a domain configuration with conduit interfacing porous media. We first recall the CHSD system derived in \cite{HSW2014}. Let $\Omega_c$ denote the conduit region and $\Omega_m$ denote the porous media. The interface between the two parts (i.e., $\partial \Omega_c\cap \partial \Omega_m$) is denoted by $\Gamma_{cm}$, on which $\mathbf{n}_{cm}$ is the unit normal to $\Gamma_{cm}$ pointing from $\Omega_c$ to $\Omega_m$. Then we define $\Gamma_c=\partial \Omega_c\backslash \Gamma_{cm}$ and $\Gamma_m=\partial \Omega_m\backslash \Gamma_{cm}$, with $\mathbf{n}_c, \mathbf{n}_m$ being the unit outer normals to $\Gamma_{c}$ and $\Gamma_{m}$. On the interface $\Gamma_{cm}$, we denote by $\{\btau_i\}$ $(i=1,...,d-1)$ a local orthonormal basis for the tangent plane to $\Gamma_{cm}$.  A two dimensional geometry is illustrated in Figure \ref{domain}.
In turn, the CHSD system takes the following form
\begin{eqnarray}
&&\rho_0 \partial_t\mathbf{u}_c = \nabla \cdot \mathbb{T}(\mathbf{u}_c, P_c) - \varphi_c\nabla\mu_c,\ \ \mbox{in}\ \Omega_c,\label{pdeSys1}\\
&&\nabla \cdot \mathbf{u}_c = 0,\ \ \mbox{in}\ \Omega_c, \\
&&\partial_t\varphi_c+\nabla\cdot(\mathbf{u}_c\varphi_c) = {\rm div}({\rm M}(\varphi_c)\nabla\mu_c),\ \ \mbox{in}\ \Omega_c,\\
&&\frac{\rho_0}{\chi}\partial_t\mathbf{u}_m + \nu(\varphi_m)\Pi^{-1}\mathbf{u}_m = -(\nabla P_m + \varphi_m\nabla\mu_m),\ \ \mbox{in}\ \Omega_m\\
&&\nabla \cdot \mathbf{u}_m = 0,\ \ \mbox{in}\ \Omega_m, \\
&&\partial_t\varphi_m + \nabla\cdot(\mathbf{u}_m\varphi_m) = {\rm div}({\rm M}(\varphi_m)\nabla\mu_m),\ \ \mbox{in}\ \Omega_m.
\end{eqnarray}
The chemical potentials $\mu_c, \mu_m$ turn out to be 
\begin{equation}\label{chem}
\mu_j = \gamma[\frac{1}{\epsilon}(\varphi_j^3-\varphi_j)-\epsilon\Delta\varphi_j],~~~j\in\{c,m\},
\end{equation}
and the Cauchy stress tensor $\mathbb{T}$ is given by
\begin{equation}
\mathbb{T}(\mathbf{u}_c,P_c) = 2\nu(\varphi_c)\mathbb{D}(\mathbf{u}_c)-P_c\mathbb{I},
\end{equation}
in which $\mathbb{D}(\mathbf{u}_c) = \frac{1}{2}(\nabla\mathbf{u}_c+\nabla\mathbf{u}_c^T)$ and $\mathbb{I}$ is the $d\times d$ identity matrix. Here $\rho_0$ is the density of the fluid, M is the mobility satisfying $0 < M_0 \leq {\rm M} \leq M_1$, $\chi$ is the porosity, $\nu$ is the viscosity satisfying $0<\nu_0\leq\nu\leq\nu_1$. In addition, we assume that  both the mobility ${\rm M}$ and the viscosity $\nu$ are Lipschitz continuous. $\Pi$ is the permeability matrix of size $d\times d$ which is assumed to be bounded, symmetric and uniformly positive definite.  The parameter $\gamma$ in \eqref{chem} is a positive constant related to the surface tension.

\ignore{\color{red}
We assume that both the mobility ${\rm M}$ and the viscosity $\nu$ satisfy Lipschitz continuity. That is, there exists a constant $C>0$ such that for all $x_1, x_2 \in \mathbb{R}$,
\begin{eqnarray}
\left|{\rm M}(x_1) - {\rm M}(x_2)\right| &\leq& C\left|x_1-x_2\right|, \\
\left|\nu(x_1) - \nu(x_2)\right| &\leq& C\left|x_1-x_2\right|.
\end{eqnarray}}
The CHSD system is subject to the following boundary and interface conditions:

\noindent\textbf{Boundary conditions on $\Gamma_c$ and $\Gamma_m$}:
\begin{eqnarray}
&&\mathbf{u}_c = \mathbf{0},\ \ \ \pp{\varphi_c}{\mathbf{n}_c} = \pp{\mu_c}{\mathbf{n}_c} = 0,\ \ \ \mbox{on}\ \Gamma_c,\\
&&\mathbf{u}_m\cdot\mathbf{n}_m = 0,\ \ \pp{\varphi_m}{\mathbf{n}_m} = \pp{\mu_m}{\mathbf{n}_m} = 0,\ \ \ \mbox{on}\ \Gamma_m.
\end{eqnarray}

\noindent\textbf{Interface conditions on $\Gamma_{cm}$}:
\begin{eqnarray}
&&\varphi_m = \varphi_c,~~\pp{\varphi_m}{\mathbf{n}_{cm}} = \pp{\varphi_c}{\mathbf{n}_{cm}},\ \ \ \mbox{on}\ \Gamma_{cm},\\
&&\mu_m = \mu_c,\ \ {\rm M}(\varphi_m)\pp{\mu_m}{\mathbf{n}_{cm}} = {\rm M}(\varphi_c)\pp{\mu_c}{\mathbf{n}_{cm}},\ \ \ \mbox{on}\ \Gamma_{cm},\\
&&\mathbf{u}_m\cdot\mathbf{n}_{cm} = \mathbf{u}_c\cdot\mathbf{n}_{cm},\ \ \ \mbox{on}\ \Gamma_{cm},\\
&&-2\nu(\varphi_c)\mathbf{n}_{cm}\cdot\mathbb{D}(\mathbf{u}_c)\mathbf{n}_{cm} + P_c = P_m,\ \ \ \mbox{on}\ \Gamma_{cm},\label{IntfPcPm}\\
&&-\nu(\varphi_c)\btau_i\cdot\mathbb{D}(\mathbf{u}_c)\mathbf{n}_{cm}= \alpha_{BJSJ}\frac{\nu(\varphi_m)}{2\sqrt{{\rm tr}(\Pi)}}\btau_i\cdot\mathbf{u}_c,\ \ \ i=1,\ldots,d-1,\mbox{on}\ \Gamma_{cm},\label{LastIntf}
\end{eqnarray}
where $\alpha_{BJSJ}$ is an empirical parameter in the Beavers-Joseph-Saffman-Jones(BJSJ) condition and ${\rm tr}(\Pi)$ is the trace of $\Pi$.

 Define the total energy of the system as follows:
\begin{equation}\label{totenergy}
\mathcal{E}(t):=\int_{\Omega_c}\frac{\rho_0}{2}|\mathbf{u}_c|^2dx + \int_{\Omega_m}\frac{\rho_0}{2\chi}|\mathbf{u}_m|^2dx + \gamma\int_{\Omega}[\frac{\epsilon}{2}|\nabla\varphi|^2+\frac{1}{\epsilon}F(\varphi)]dx,
\end{equation}
where $F(\varphi)=\frac{1}{4}(\varphi^2-1)^2$.  The CHSD system \eqref{pdeSys1}-\eqref{LastIntf} obeys a dissipative energy law \cite{CHW2017}:
\be
\frac{d}{dt} \mathcal{E}(t)=-\mathcal{D}(t) \le 0,\quad \forall\, t\geq 0, \label{EnergyLaw}
\ee
 where the rate of energy dissipation $\mathcal{D}$ is given by
\bea
\mathcal{D}(t)&=& \int_{\Omega_m} \nu(\varphi_m)\Pi^{-1}|\ub_m|^2dx
+\int_{\Omega_c}2\nu(\vp_c)|\mathbb{D}(\ub_c)|^2dx
\non\\
&& +\int_{\Omega}{\rm M}(\vp)|\nabla\mu(\vp)|^2dx
+\int_{\Gamma_{cm}}\alpha_{BJSJ}\frac{ \nu(\vp)}{\sqrt{{\rm trace}(\Pi)}}\sum_{i=1}^{d-1}|\ub_c\cdot\btau_i|^2 dS \geq 0 . \label{D}
\eea

The CHSD system \eqref{pdeSys1}-\eqref{LastIntf} is systematically derived via Onsager's extremum principle in \cite{HSW2014}. Well-posedness of a variant of the CHSD model is studied in \cite{HWW2014}. A decoupled unconditionally stable numerical algorithm for solving the CHSD system is proposed in \cite{CHW2017}. Here we focus on the error analysis of a similar decoupled numerical scheme (cf. Sec. \ref{NuSch}) in which the computation of Stokes equations and Darcy equations are nevertheless coupled. The decoupling between the Cahn-Hilliard equation and fluid equations is accomplished by a special technique of operator splitting in which an intermediate velocity for advection in the Cahn-Hilliard equation is defined in terms of the capillarity from fluid equations. Application of this specific fractional step method for solving  phase field models is first reported in \cite{Minjeaud2013},  and later in \cite{ShYa2015}. To the best of our knowledge, error analysis of the decoupled scheme via the aforementioned operator splitting has not been reported elsewhere for any phase field model coupled with fluid motion.

There have been some convergence analysis works for either the Cahn-Hilliard-Navier-Stokes (Stokes) (CHNS, CHS) or the Cahn-Hilliard-Darcy (Hele-Shaw) system (CHD, CHHS) in recent years. The convergence of certain finite element numerical solutions to weak solutions of the CHNS equations was proved in \cite{Feng2006}, and a similar analysis is perform for the CHHS system in \cite{FeWi2012}. In \cite{DFW2015} the authors have established optimal convergence rates for a mixed finite element method for solving the CHS system, with first order temporal accuracy. More recently, an optimal rate error estimate is presented for a second-order accurate numerical scheme for solving the CHNS equations in \cite{DWWW2017}.  A similar error estimate was also reported in \cite{CaSh2018}, based on a finite element discretization of a linear, weakly coupled energy stable scheme  for the CHNS system.  As for the CHHS system, in which the kinematic diffusion term is replaced by a damping one, optimal error analysis has been presented in \cite{CLWW2016, LCWW2017}, in the framework of finite difference and finite element spatial approximations, respectively.

The CHSD system consist of the CHS and the CHD equations, coupled together via a set of domain interface boundary conditions. Hence the advection in the Cahn-Hilliard flow is involved with both the Stokes and the Darcy velocity fileds. While the Stokes velocity has a regularity of $L^2(0, T; H^1)$, the Darcy velocity is only of $L^\infty(0, T; L^2)$. With the $L^2 (0, T; H^1)$ bound of the velocity field, a uniform maximum norm estimate of the phase has been derived, which significantly simplifies the error analysis for the CHNS system \cite{DWWW2017} and the CHS equations \cite{DFW2015}. On the other hand, for the CHD system, only an  $L^p (0, T; L^\infty)$ bound (with a finite value of $p$) could be established for the phase variable, as analyzed in~\cite{LCWW2017}. The lack of uniform bound of the phase variable has dramatically complicated the error analysis of the nonlinear advection associated with the Cahn-Hilliard equation. A similar difficulty is encountered here for the error analysis of the CHSD system. To overcome this subtle difficulty, we perform an $L^2 (0,T; H^3)$ bound estimate of the phase variable in the numerical solution, which is accomplished by the usage of a discrete Gagliardo-Nirenberg inequality in the finite element setting. This bound will play an important role to pass through the error estimate. Such a technique has been applied in the analysis for the CHHS system in the existing literature, as reported in~\cite{CLWW2016, LCWW2017, chen19a}. Moreover, the CHSD system contains a coupling between the CHS and CHD equations, the corresponding estimates are expected to be even more challenging than the ones for the CHHS model. 

The rest of the article is organized as follows. In Section~\ref{NuSch} we introduce the weak formulation of the CHSD system and present the decoupled numerical scheme. Some preliminary analysis including the stability estimates are gathered in Section~\ref{Prel}. The detailed error analysis of the numerical scheme is carried out in Section~\ref{Err-ana}. Finally, some concluding remarks are provided in Section~\ref{sec:conclusion}.

\section{The numerical scheme}\label{NuSch}

\subsection{The weak formulation }\label{WF}

For the CHSD problem we introduce the following spaces
\bea
 \BH({\rm div}; \Omega_j) &:=&\{\mathbf{w}\in \mathbf{L}^2(\Omega_j)~|~\nabla \cdot \mathbf{w}\in \mathbf{L}^2(\Omega_j)\}, \quad \quad j\in \{c,m\},\non\\
\mathbf{H}_{c,0}&:=&\{\mathbf{w}\in \BH^1(\Omega_c)~|~\mathbf{w}=\mathbf{0}\text{ on
}\Gamma_{c}\},
 \non\\
\mathbf{H}_{c,\text{div}}&:=&\{
\mathbf{w}\in\mathbf{H}_{c,0}~|~\nabla \cdot\mathbf{w}=0\}, \nonumber
 \\
\mathbf{H}_{m,0}&:=&\{\mathbf{w}\in \BH({\rm div}; \Omega_m)~|~\mathbf{w}\cdot \mathbf{n}_m=0 \ \text{on}\ \Gamma_{m}\}, \non\\
\mathbf{H}_{m,\mathrm{div}}&:=&\{\mathbf{w}\in\mathbf{H}_{m,0}~|~\nabla \cdot\mathbf{w}=0\},\non\\
 X_m&:=& H^1(\Omega_m) \cap L^2_0(\Omega_m).\non
 \eea
 Here $L^2_0(\Omega_m)$ is a subspace of $L^2$ whose elements are of mean zero. We also use the notation $L^2_0(\Omega)$ which is defined similarly and will be used later. We denote $(\cdot, \cdot)_c$, $(\cdot, \cdot)_m$ the inner products on the spaces $L^2(\Omega_c)$, $L^2(\Omega_m)$, respectively (also for the corresponding vector spaces). The inner product on $L^2(\Omega)$ is simply denoted by $(\cdot, \cdot)$. In turn, it is clear that
$$
(u,v)=(u_m,v_m)_m+(u_c,v_c)_c, \quad \|u\|_{L^2(\Omega)}^2=\|u_m\|_{L^2(\Omega_m)}^2+\|u_c\|_{L^2(\Omega_c)}^2,
$$
where $u_m:=u|_{\Omega_m}$ and $u_c:=u|_{\Omega_c}$. We will suppress the dependence on the domain in the $L^2$ norm if there is no ambiguity. And also, $H'$ stands for the dual space of $H$ with the duality induced by the $L^2$ inner product. For simplicity, we denote $\|\cdot\|:=\|\cdot\|_{L^2}$, and $\|\cdot\|_p:=\|\cdot\|_{L^p}$ for $1\leq p\leq \infty,\ p\neq 2$. In addition, the notation $\left\|\cdot\right\|_{cm}$ is introduced as the $L^2$ norm on the interface $\Gamma_{cm}$. For all the functions $f$, $\overline{f}$ represents the mean value of $f$ on its domain.

The definition of the weak formulation of the 3-D CHSD system is given below. The 2-D case could be similarly defined with slight changes in time integrability of the functions.
\begin{mydef}\label{defweak}
Suppose that $d=3$ and $T>0$ is arbitrary.
We consider the initial data $\vp_{0}\in H^1(\Omega), \ub_c(0) \in \mathbf{H}_{c,{\rm div}}, \ub_m(0) \in \mathbf{H}_{m,{\rm div}}$. The functions $(\ub_c, P_c, \mathbf{u}_m, P_m, \vp, \mu)$ with the following properties
\bea
&  \ub_c\in  L^\infty(0, T; \mathbf{L^2}(\Omega_c)) \cap L^2(0, T; \mathbf{H}_{c,0}),  \frac{\partial \ub_c}{\partial t} \in L^{\frac{4}{3}}(0,T; (\mathbf{H}_{c,0})^\prime),  \\
& \ub_m\in L^\infty(0, T; \mathbf{L^2}(\Omega_m)) \cap L^2(0, T; \mathbf{H}_{m,0}),\frac{\partial \ub_m}{\partial t} \in L^{\frac{4}{3}}(0,T; (\mathbf{H}_{m,0})^\prime), \label{regpm}\\
& P_c \in L^\frac43(0,T; L^2(\Omega_c)), \quad P_m \in L^\frac43(0,T; X_m), \\
& \vp\in L^\infty(0,T; H^1(\Omega))\cap L^2(0,T; H^3(\Omega)), \vp_t \in  L^2(0;T; (H^1(\Omega))'),\\
& \mu\in L^2(0, T; H^1(\Omega)),
\eea
is called a finite energy weak solution of the CHSD system \eqref{pdeSys1}--\eqref{LastIntf}, if the following conditions are satisfied:

(1) For any $v, \phi\in  H^1(\Omega)$,
\begin{eqnarray}
&&\langle \partial_t\varphi,v)+({\rm M}(\vp) \nabla \mu(\vp),\nabla v)-(\ub \varphi ,\nabla v)=0, \label{weak_CH1}\\
&&\gamma \left[\frac{1}{\epsilon}(f(\varphi),\phi)+\epsilon(\nabla\varphi,\nabla \phi) \right]-(\mu(\varphi),\phi)=0,  \quad  f(\varphi):= \varphi^3 - \varphi .  \label{weak_CH2}
\end{eqnarray} 

(2) For any $\mathbf{v}_c\in  \mathbf{H}_{c,0}$ and $q_c\in  L^2(\Omega_c)$,
 \begin{eqnarray}
 &&\rho_0\langle \partial_t \ub_c,\mathbf{v}_c \rangle_c + a_c(\mathbf{u}_{c},\mathbf{v}_{c}) + b_c(\mathbf{v}_{c},P_{c}) + \int_{\Gamma_{cm}} P_m (\mathbf{v}_c\cdot \mathbf{n}_{cm}) dS\non\\
&&  - b_c(\ub_{c},q_{c}) +(\vp_c \nabla \mu(\varphi_c),  \mathbf{v}_c )_c =0,
 \label{weak_S}
\end{eqnarray}
where
\begin{eqnarray*}
a_c(\mathbf{u}_{c},\mathbf{v}_{c}) = 2\left(\nu(\vp_c)\mathbb{D}(\mathbf{u}_{c}),\mathbb{D}(\mathbf{v}_{c})\right)_c &+& \sum_{i=1}^{d-1}\int_{\Gamma_{cm}}\alpha_{BJSJ}\frac{\nu(\vp)}{\sqrt{{\rm tr} (\Pi)}}
(\mathbf{u}_{c}\cdot\btau_i)(\mathbf{v}_{c}\cdot\btau_i)dS\ \ , \\
b_c(\mathbf{v}_{c},q_{c}) &=& -(\nabla\cdot\mathbf{v}_{c},q_{c})_c.
\end{eqnarray*}

(3) For any $\mathbf{v}_m\in  \mathbf{H}_{m,0}$ and $q_m\in  H^1(\Omega_m)$,
\begin{eqnarray}
&&\frac{\rho_0}{\chi} \langle \partial_t \ub_m,\mathbf{v}_m \rangle_m + a_m(\mathbf{u}_{m},\mathbf{v}_{m}) + b_m(\mathbf{v}_{m}, P_{m}) - b_m(\ub_{m}, q_{m})\non \\
&& +(\vp_m \nabla \mu(\varphi_m),  \mathbf{v}_m )_m -\int_{\Gamma_{cm}} \ub_c \cdot \mathbf{n}_{cm} q_m\, ds=0,
\label{weak_D}
\end{eqnarray}
where
\begin{eqnarray*}
a_m(\mathbf{u}_{m},\mathbf{v}_{m})&=&\big(\nu(\vp_m)\Pi^{-1}\mathbf{u}_{m}, \mathbf{v}_{m}\big)_m\ \ ,\\
b_m(\mathbf{v}_{m}, q_{m})&=&( \mathbf{v}_{m}, \nabla q_{m})_m\ \ .
\end{eqnarray*}

(4)  $\vp|_{t=0}=\vp_{0}(x), \ub_c|_{t=0}=\ub_c(0), \ub_m|_{t=0}=\ub_m(0).$

(5) The finite energy solution satisfies the energy inequality
\be
\mathcal{E}(t)  +\int_s^t\mathcal{D}(\tau) d\tau \leq \mathcal{E}(s), \label{Energyinq}
\ee
for all $t\in [s,T)$ and almost all $s\in [0,T)$ (including $s=0$), where the total energy $\mathcal{E}$ is given by \eqref{totenergy}.

\end{mydef}
\subsection{The numerical scheme}
Let $\tau>0$ be the time step size,  $K=[T/\tau]$, and set $t^{k}=k\tau$ for $0\leq k \leq K$. Similarly, we denote $\ub^k$ as a numerical approximation to $\ub(t^k)=\ub(k\tau)$, with a notation $\ub(t):=\ub(\cdot,t)$ for simplicity. Let $\mathcal{T}_c^h$ and $\mathcal{T}_m^h$ be a quasi-uniform triangulation of the domain $\Omega_c$ and $\Omega_m$ with mesh size $h$. Then $\mathcal{T}^h:=\mathcal{T}_c^h\cup \mathcal{T}_m^h$ forms a triangulation of the whole domain $\Omega$. $\mathcal{T}_c^h$ and $(\mathcal{T}_m^h)$ coincide on the interface $\Gamma_{cm}$. Let $Y_h$ denote the finite element approximation of $H^1(\Omega)$, such as
\begin{eqnarray}
Y_h=\{v_h\in C(\bar{\Omega})\big|v_h|_K \in P_r(K),\forall K\in\mathcal{T}_h\}\nonumber.
\end{eqnarray}
Additionally, we introduce $\mathring{Y}_h:=Y_h\cap L_0^2(\Omega)$. Let $\mathbf{X}_c^h, M_c^h, \mathbf{X}_m^h, M_m^h$ be the finite element approximation of $\mathbf{H}_{c,0},L^2(\Omega_c),\mathbf{H}_{m,0},X_m$ respectively, while the approximation polynomials have adequate degrees. We assume that $\mathbf{X}_c^h$ and $M_c^h$ are stable approximation spaces for Stokes velocity and pressure in the sense that
\begin{align}\label{inf-supS}
\sup_{\mathbf{v}_h \in \mathbf{X}_c^h} \frac{(\nabla \cdot \mathbf{v}_h,  q_h)_c}{||\mathbf{v}_h||_{H^1}} \geq c||q_h||, \quad \forall q_h \in M_c^h.
\end{align}
The validity of such an inf-sup condition for some standard finite element spaces can be found in \cite{layton2002coupling}. The classical P2-P0,   Taylor-Hood finite element spaces and the Mini finite element spaces    are commonly adopted in practice for $\mathbf{X}_c^h$ and $M_c^h$, cf. \cite{layton2002coupling}, \cite{GiRa1986}. The spaces $\mathbf{X}_m^h$ and $M_m^h$ are assumed to be stable in the sense that 
\begin{align}\label{inf-supD}
\sup_{\mathbf{v}_h \in \mathbf{X}_m^h} \frac{( \mathbf{v}_h, \nabla q_h)_m}{||\mathbf{v}_h||} \geq c||q_h||, \quad \forall q_h \in M_m^h.
\end{align}
In particular, we notice that the Taylor-Hood finite element spaces satisfy the above condition.

We will focus on the error analysis of the following unconditionally energy stable scheme that decouples the computation of the Cahn-Hilliard flow from that of fluid equations, i.e.,  for a totally decoupled scheme; see the related descriptions in~\cite{CHW2017}. FGiven $0\leq k\leq K-1$, find $\big(\varphi_h^{k+1}, \mu_h^{k+1}, \ub_{c, h}^{k+1}, P_{c,h}^{k+1}, \ub_{m, h}^{k+1}, P_{m,h}^{k+1}\big)$ $\in$ $Y_h \times Y_h \times \mathbf{X}_c^h \times M_c^h \times \mathbf{X}_m^h \times M_m^h$ such that for all $\left(v, \phi, \mathbf{v}_{c}, q_{c}, \mathbf{v}_{m}, q_{m}\right)$ $\in$ $Y_h \times Y_h \times \mathbf{X}_c^h \times M_c^h \times \mathbf{X}_m^h \times M_m^h$ there holds
\begin{subequations}
	\begin{eqnarray}
	&&( \delta_t\varphi^{k+1}_h,v)+({\rm M}(\vp^k_h) \nabla \mu^{k+1}_h,\nabla v)-(\baru^{k+1}_h \varphi^k_h ,\nabla v)=0, \label{PD_CH1}\\
	&&\gamma \left[\frac{1}{\epsilon}(f(\varphi^{k+1}_h, \varphi^k_h),\phi)+\epsilon(\nabla\varphi^{k+1}_h,\nabla \phi) \right]-(\mu^{k+1}_h,\phi)=0, \label{PD_CH2} \\
	&&\rho_0( \delta_t\ub_{c,h}^{k+1},\mathbf{v}_{c} )_c  +a_c^k(\ub_{c,h}^{k+1},\mathbf{v}_{c}\big)+ b_c(\mathbf{v}_{c}, P_{c,h}^{k+1})
	+  \int_{\Gamma_{cm}} P_{m,h}^{k+1} (\mathbf{v}_{c}\cdot \mathbf{n}_{cm}) dS \non\\
	&&\ \ \ \  - b_c(\mathbf{u}^{k+1}_{c,h}, q_{c}) +(\vp_{c,h}^k \nabla \mu_{c,h}^{k+1},  \mathbf{v}_{c} )_c =0,
	\label{PD_S} \\
	&&\frac{\rho_0}{\chi} (\delta_t\ub_{m,h}^{k+1},\mathbf{v}_{m} )_m + a_m^k\big(\mathbf{u}_{m,h}^{k+1},\mathbf{v}_{m}\big)
	+b_m(\mathbf{v}_{m}, P_{m,h}^{k+1})  +(\vp_{m,h}^k \nabla \mu_{m,h}^{k+1},  \mathbf{v}_{m} )_m \non \\
	&&\ \ \ \ -\int_{\Gamma_{cm}} \ub_{c,h}^{k+1} \cdot \mathbf{n}_{cm} q_{m}\, ds- b_m(\ub_{m,h}^{k+1},  q_{m})=0,
	\label{PD_D}
	\end{eqnarray}
\end{subequations}
where
\begin{eqnarray}
&& f(\varphi^{k+1}_h, \varphi^k_h):=(\varphi^{k+1}_h)^3-\varphi^k_h, \ \ \ \ \delta_t\varphi^{k+1}_h:= \frac{\varphi^{k+1}_h-\varphi^{k}_h}{\tau}, \\
&& \baru_h ^{k+1}=\begin{cases}
&\baru_{m,h}^{k+1}, \quad x \in \Omega_m, \\
&\baru_{c,h}^{k+1}, \quad x \in \Omega_c,
\end{cases},\ \ \ \
\begin{cases}
& \frac{\rho_0}{\chi}\frac{\baru_{m,h}^{k+1}-\ub_{m,h}^k}{\tau}+\varphi_{m,h}^{k}\nabla \mu_{m,h}^{k+1}=0, \\
& \rho_0 \frac{\baru_{c,h}^{k+1}-\ub_{c,h}^k}{\tau}+\varphi_{c,h}^{k}\nabla \mu_{c,h}^{k+1}=0,
\end{cases} \label{InterM} ,\\
&& a_c^k(\ub_{c,h}^{k+1},\mathbf{v}_{c})=2(\nu(\vp_{c,h}^k)\mathbb{D}(\ub_{c,h}^{k+1}),\mathbb{D}(\mathbf{v}_{c}))_c \non\\
&&\ \ \ \ \ \ \ \  +\sum_{i=1}^{d-1}\int_{\Gamma_{cm}}\alpha_{BJSJ}\frac{\nu(\vp_{c,h}^k)}{\sqrt{{\rm tr} (\Pi)}}(\ub_{c,h}^{k+1}\cdot\btau_i)(\mathbf{v}_{c}\cdot\btau_i)dS, \label{AC}\\
&& b_c(\mathbf{v}_{c}, q_{c})=-(\nabla \cdot \mathbf{v}_{c}, q_{c})_c, \label{BC}\\
&& a_m^k(\mathbf{u}_{m,h}^{k+1},\mathbf{v}_{m})=\big(\nu(\varphi_{m,h}^{k})\Pi^{-1}\mathbf{u}_{m,h}^{k+1}, \mathbf{v}_{m}\big)_m, \label{AM}\\
&& b_m(\mathbf{v}_{m}, q_{m})=( \mathbf{v}_{m}, \nabla q_{m})_m. \label{BM}
\end{eqnarray}
The initial values are taken as follows: 
\begin{eqnarray}\label{initialC}
\varphi_h^0 = \mathcal{P} \varphi^0,\ \ \ \ub_{j,h}^0 = \mathcal{P}_{j,u}^{0} \ub_j^0,\ \ j\in\{c,m\}.
\end{eqnarray}

The unique solvability of the proposed scheme \eqref{PD_CH1}--\eqref{BM} has been  proved via a  convexity analysis, and the energy stability is ensured by a careful estimate; the details could be found in~\cite{CHW2017}. In this article, we focus on the optimal rate convergence analysis and error estimate. 

\section{Some preliminary estimates}\label{Prel}
Some projections are needed in the later analysis: \\
Ritz projection $\mathcal{P}: H^{1}(\Omega)\rightarrow Y_h$,
\begin{eqnarray}
\big(\nabla(\mathcal{P}\varphi-\varphi),\nabla v\big)=0,\ \ \forall v\in Y_h,\ \ \ (\mathcal{P}\varphi-\varphi,1)=0,
\end{eqnarray}
and for $\phi=\vp(t), \forall t\in[0,T]$, where $\vp$ is of the weak solution to CHSD system \eqref{pdeSys1}--\eqref{LastIntf}, we define the modified Ritz projection $\widetilde{\mathcal{P}}^{\phi}$: $H^{1}(\Omega)\rightarrow Y_h$,
\begin{eqnarray}
\big({\rm M}(\phi)\nabla(\widetilde{\mathcal{P}}^{\phi}\mu-\mu),\nabla v\big)=0,\ \ \forall v\in Y_h,\ \ \ (\widetilde{\mathcal{P}}^{\phi}\mu-\mu,1)=0 . 
\end{eqnarray}
Stokes--Darcy projection $\left(\mathcal{P}_{c,u}^{\phi},\mathcal{P}_{c,p}^{\phi},\mathcal{P}_{m,u}^{\phi},\mathcal{P}_{m,p}^{\phi}\right)$: $\left(\mathbf{H}_{c,0},L^2(\Omega_c),\mathbf{H}_{m,0},X_m\right)\rightarrow \big(\mathbf{X}_c^h, M_c^h,$ $\mathbf{X}_m^h, M_m^h\big)$, which, for all $\mathbf{v}_c \in\mathbf{X}_c^h, q_c \in M_c^h, \mathbf{v}_m \in\mathbf{X}_m^h, q_m \in M_m^h$, satisfies the following equalities:
\begin{eqnarray}
&& 2\Big(\nu(\phi_c)\mathbb{D}\big(\mathcal{P}_{c,u}^{\phi}\ub_c\big),\mathbb{D}(\mathbf{v}_c)\Big)_c + \sum_{i=1}^{d-1}\int_{\Gamma_{cm}}\alpha_{BJSJ}\frac{\nu(\phi_c)}{\sqrt{{\rm tr}(\Pi)}} \left((\mathcal{P}_{c,u}^{\phi}\ub_c)\cdot\btau_i\right)\left(\mathbf{v}_c\cdot\btau_i\right)dS \nonumber\\
&&\ \ \ \ \ - \left(\mathcal{P}_{c,p}^{\phi}P_c,\nabla\cdot\mathbf{v}_c\right)_c + \int_{\Gamma_{cm}}(\mathcal{P}_{m,p}^{\phi}P_m)\left(\mathbf{v}_c\cdot\mathbf{n}_{cm}\right)dS + \left(\nabla\cdot\big(\mathcal{P}_{c,u}^{\phi}\ub_c\big),q_c\right)_c \nonumber\\
&=& 2\Big(\nu(\phi_c)\mathbb{D}(\ub_c),\mathbb{D}(\mathbf{v}_c)\Big)_c + \sum_{i=1}^{d-1}\int_{\Gamma_{cm}}\alpha_{BJSJ}\frac{\nu(\phi_m)}{\sqrt{{\rm tr}(\Pi)}} \left(\ub_c\cdot\btau_i\right)\left(\mathbf{v}_c\cdot\btau_i\right)dS \nonumber\\
&&\ \ \ \ \ - \left(P_c,\nabla\cdot\mathbf{v}_c\right)_c + \int_{\Gamma_{cm}}P_m\left(\mathbf{v}_c\cdot\mathbf{n}_{cm}\right)dS + \left(\nabla\cdot\ub_c,q_c\right)_c\ \ ,
\end{eqnarray}
\begin{eqnarray}
&& \Big(\nu(\phi_m)\Pi^{-1}\big(\mathcal{P}_{m,u}^{\phi}\ub_m\big),\mathbf{v}_m\Big)_m + \left(\nabla\big(\mathcal{P}_{m,p}^{\phi}P_m\big),\mathbf{v}_m\right)_m - \Big(\mathcal{P}_{m,u}^{\phi}\ub_m,\nabla q_m\Big)_m - \int_{\Gamma_{cm}}\big(\mathcal{P}_{c,u}^{\phi}\ub_c\big)\cdot\mathbf{n}_{cm}q_m dS \nonumber\\
&=& \Big(\nu(\phi_m)\Pi^{-1}\ub_m,\mathbf{v}_m\Big)_m + \left(\nabla P_m,\mathbf{v}_m\right)_m - \left(\ub_m,\nabla q_m\right)_m - \int_{\Gamma_{cm}}\ub_c\cdot\mathbf{n}_{cm}q_m dS.
\end{eqnarray}
Especially, for $0\leq k\leq K$, we rewrite the notation of the projections above as follows:
\begin{eqnarray}
\widetilde{\mathcal{P}}^k &:=& \widetilde{\mathcal{P}}^{(\vp^k)}, \\
\left(\mathcal{P}_{c,u}^{k},\mathcal{P}_{c,p}^{k},\mathcal{P}_{m,u}^{k},\mathcal{P}_{m,p}^{k}\right) &:=& \left(\mathcal{P}_{c,u}^{(\vp^k)},\mathcal{P}_{c,p}^{(\vp^k)},\mathcal{P}_{m,u}^{(\vp^k)},\mathcal{P}_{m,p}^{(\vp^k)}\right).
\end{eqnarray}

What follows is a standard result of Ritz projection \cite{BrSc2008}. There exists a constant $C>0$ depending on $M_0, M_1$, such that the Ritz projections $\mathcal{P}$ and $\widetilde{\mathcal{P}}^k$ satisfies
\begin{eqnarray}
&& \big\|\mathcal{P}\varphi-\varphi\big\|_p + h\big\|\nabla(\mathcal{P}\varphi-\varphi)\big\|_{p} \leq C h^{q+1}\big\|\varphi\big\|_{W_p^{q+1}}, \\
&& \big\|\widetilde{\mathcal{P}}^k\varphi-\varphi\big\| + h\big\|\nabla(\widetilde{\mathcal{P}}^k\varphi-\varphi)\big\| \leq C h^{q+1}\big\|\varphi\big\|_{H^{q+1}},
\end{eqnarray}
for all $\varphi \in H^{q+1}(\Omega),\ q\geq 0,\ p\in[2,\infty]$, and all $0\leq k\leq K$ with $Y_h$ consisting of polynomials of order q or higher.

For the Stokes-Darcy projection, the following error estimates have been established in \cite{CGSW2013, MuZh2010, RiYo2005}
\begin{eqnarray}
\big\|\mathbf{u}_c-\mathcal{P}_{c,u}^k\mathbf{u}_c\big\|_{H^1(\Omega_c)}+ \big\|\mathbf{u}_m-\mathcal{P}_{m,u}^k\mathbf{u}_m\big\| \leq h^q \left(\big\| \mathbf{u}_c\big\|_{H^{q+1}(\Omega_c)}+\big\| \mathbf{u}_m\big\|_{H^{q+1}(\Omega_m)}\right).
\end{eqnarray}

Here we introduce the linear operator $\textsf{T}_h: \mathring{Y}_h\rightarrow \mathring{Y}_h$, which is defined via the variational problem: given $\zeta\in\mathring{Y}_h$, find $\textsf{T}_h(\zeta)\in\mathring{Y}_h$ such that
\be
\big(\nabla\textsf{T}_h(\zeta), \nabla\xi\big)=\big(\zeta, \xi\big),\ \ \ \forall \xi\in\mathring{Y}_h.
\ee
With this operator, we are able to define the following $\left\|\cdot\right\|_{-1,h}$ norm:
\be
\left\|\zeta\right\|_{-1,h} := \left\|\nabla \textsf{T}_h(\zeta)\right\| = \sqrt{\big(\nabla \textsf{T}_h(\zeta), \nabla \textsf{T}_h(\zeta)\big)} = \sqrt{\big(\zeta, \textsf{T}_h(\zeta)\big)},\ \ \ \forall \zeta\in\mathring{Y}_h.
\ee

We also define the discrete Laplacian, $\Delta_h$: $Y_h\rightarrow\mathring{Y}_h$ as follows:
for any $v_h\in Y_h,\ \Delta_h v_h\in\mathring{Y}_h$ denotes the unique solution to the problem
\be
(\Delta_h v_h, \xi) = - (\nabla v_h, \nabla \xi), \ \ \ \forall\xi\in Y_h.
\ee

We recall the following discrete Gagliardo--Nirenberg inequality from \cite{HeRa1982, LCWW2017} which is needed for the uniform estimate of the order parameter $\varphi_h^{k+1}$.
\begin{lemma} \label{nirenberg}
Suppose that $\Omega$ is a convex and polyhedral domain. Then, for any $\varphi_h \in Y_h$,
\begin{eqnarray}
\big\|\varphi_h\big\|_{L^{\infty}} \leq C\big\|\Delta_h\varphi_h\big\|^{\frac{d}{2(6-d)}}\big\|\varphi_h\big\|_{L^6}^{\frac{3(4-d)}{2(6-d)}} + C\big\|\varphi_h\big\|_{L^6}\ ,\ \ \forall \varphi_h \in Y_h.\ \ \ d=2,3,
\end{eqnarray}
and consequently,
\begin{eqnarray}
\left\|\vp_h-\overline{\vp_h}\right\|_{L^{\infty}} \leq C \left\|\nabla\Delta_h\vp_h\right\|^{\frac{d}{4(6-d)}} \left\|\nabla\vp_h\right\|^{\frac{24-5d}{4(6-d)}} + C\left\|\nabla\vp_h\right\|,\ \ \ d=2,3,
\end{eqnarray}
where $\overline{\vp_h}$ is the mean value of $\vp_h$.
\end{lemma}

The following technical lemma has been proved in \cite{DWWW2017}.
\begin{lemma}\label{nabla_-1,h}
Suppose $g\in H^1(\Omega)$ and $v\in \mathring{Y}_h$. Then
\begin{eqnarray}
\left|(g, v)\right|\leq C\|\nabla g\|\|v\|_{-1,h}
\end{eqnarray}
holds for some $C>0$ that is independent of $h$.
\end{lemma}
We also recall the inverse inequality
\begin{eqnarray}
\|\vp_h\|_{W^m_q} \leq C h^{d/q-d/p} h^{l-m}\|\vp_h\|_{W^l_p},\ \ \forall \vp_h\in Y_h,
\end{eqnarray}
for all $1\leq p\leq q\leq\infty,\ 0\leq l\leq m\leq1$.

The following trace theorem is necessary for the estimate of certain interface boundary terms.
\begin{lemma}\label{trace}
Suppose $\vb\in H^1(\Omega)$. Then
\begin{eqnarray}
\|\vb\|_{L^4(\partial\Omega)} \leq C\|\vb\|_{H^1(\Omega)}.
\end{eqnarray}
In particular for $\ub_h \in \mathbf{H}_{c,0}$, there holds
\begin{eqnarray}
\|\ub_h\|_{L^4(\Gamma_{cm})} \leq C \|\mathbb{D}(\ub_h)\|_{L^2(\Omega_c)}.
\end{eqnarray}
\end{lemma}

Now we derive some stability estimate of the scheme \eqref{PD_CH1}--\eqref{initialC}.  The following estimates are direct consequence of the discrete energy law established in 
 \cite{CHW2017}.
\begin{lemma}\label{energy_es}
Let $(\vp_h^{k+1}, \mu_h^{k+1}, \ub_{c,h}^{k+1}, P_{c,h}^{k+1}, \ub_{m,h}^{k+1}, P_{m,h}^{k+1})$ $\in$ $Y_h \times Y_h \times \mathbf{X}_c^h \times M_c^h \times \mathbf{X}_m^h \times M_m^h$ be the unique solution of \eqref{PD_CH1}--\eqref{initialC} for $0\leq k \leq K-1$. Then there exists a constant $C>0$ dependent on the initial data such that
\begin{eqnarray}
\max_{0\leq k\leq K} \left[\big\|\ub_{c,h}^{k}\big\|^2 + \big\|\ub_{m,h}^{k}\big\|^2 + \big\|(\varphi_h^{k})^2-1\big\|^2 + \big\|\nabla\varphi_h^{k}\big\|^2\right] &\leq& C, \\
\max_{0\leq k\leq K} \left\|\varphi_h^{k}\right\|_{H_1} &\leq& C, \label{varphi_h_H1}\\
\sum_{k=0}^{K-1} \Big[\tau\big\|\nabla\mu_h^{k+1}\big\|^2 + \tau a_c^k(\ub_{c,h}^{k+1},\ub_{c,h}^{k+1}) + \tau\big\|\ub_{m,h}^{k+1}\big\|^2 + \big\|\ub_{m,h}^{k+1}-\ub_{m,h}^{k}\big\|^2 && \nonumber\\
\ \ \ \ \ \ \ + \big\|\ub_{c,h}^{k+1}-\ub_{c,h}^{k}\big\|^2 + \big\|\nabla(\varphi_h^{k+1}-\varphi_h^{k})\big\|^2\Big] &\leq& C,  \label{energy_law}
\end{eqnarray}
hold for every $0\leq k \leq K-1$, $d=2,3$.
\end{lemma}

For the error analysis, we also need the uniform bound of the order parameter and the chemical potential for which we derive the following stability estimates,  see also  Lemma 2.13 from \cite{DFW2015}.
\begin{lemma}\label{numerical_estimate}
Let $(\vp_h^{k+1}, \mu_h^{k+1}, \ub_{c,h}^{k+1}, P_{c,h}^{k+1}, \ub_{m,h}^{k+1}, P_{m,h}^{k+1}) \in Y_h \times Y_h \times \mathbf{X}_c^h \times M_c^h \times \mathbf{X}_m^h \times M_m^h$ be the unique solution of \eqref{PD_CH1}--\eqref{initialC} for $0\leq k \leq K-1$. Then there exists some constant $C>0$ dependent on $\gamma$ and $\epsilon$ such that
\begin{eqnarray}
&& \big\|\Delta_h \varphi_h^{k+1}\big\|^2 \leq C\big\|\mu_h^{k+1}\big\|^2 + C, \label{lemma_est_1}\\
&& \big\|\mu_h^{k+1}\big\|^2 \leq \big\|\nabla\mu_h^{k+1}\big\|^2 + C, \label{lemma_est_2}\\
&& \tau \sum_{k=0}^{K-1}\left[\big\|\Delta_h \varphi_h^{k+1}\big\|^2 + \big\|\mu_h^{k+1}\big\|_{H^1}^2\right] \leq C(T+1), \label{lemma_est_3}\\
&& \tau \sum_{k=0}^{K-1} \big\|\varphi_h^{k+1}\big\|_{{\infty}}^{\frac{4(6-d)}{d}} \leq C(T+1), \label{varphi_h^4}\\
&& \tau \sum_{k=0}^{K-1} \left[\big\|\nabla\Delta_h \varphi_h^{k+1}\big\|^2 + \left\|\vp_h^{k+1}\right\|_{{\infty}}^{\frac{8(6-d)}{d}}\right] \leq C(T+1) \label{varphi_h^8},
\end{eqnarray}
hold for every $0\leq k \leq K-1$, $d=2,3$.
\end{lemma}
\begin{proof}
Setting $\phi_h=\Delta_h\vp_h^{k+1}$ in \eqref{PD_CH2}, by the uniform bound of $\left\|\vp_h^{k+1}\right\|_{H^1}$ and $\left\|\vp_h^k\right\|$ in Lemma \ref{energy_es}, we have
\begin{eqnarray}
\|\Delta_h\vp_h^{k+1}\|^2 &=& -(\nabla\vp_h^{k+1}, \nabla\Delta_h\vp_h^{k+1}) \non\\
&=& \frac{1}{\epsilon^2}\left(f(\vp_h^{k+1},\vp_h^k), \Delta_h\vp_h^{k+1}\right) - \frac{1}{\gamma\epsilon}(\mu_h^{k+1},\Delta_h\vp_h^{k+1}) \non\\
&\leq& \frac{1}{\epsilon^2}\left\|f(\vp_h^{k+1},\vp_h^k)\right\| \left\|\Delta_h\vp_h^{k+1}\right\| + \frac{1}{\gamma\epsilon}\left\|\mu_h^{k+1}\right\| \left\|\Delta_h\vp_h^{k+1}\right\| \non\\
&\leq& \frac{1}{\epsilon^2}\left(\left\|\vp_h^{k+1}\right\|_{L^6}^3 + \left\|\vp_h^k\right\|\right) \left\|\Delta_h\vp_h^{k+1}\right\| + \frac{1}{\gamma\epsilon}\left\|\mu_h^{k+1}\right\| \left\|\Delta_h\vp_h^{k+1}\right\| \non\\
&\leq& \frac{1}{\epsilon^2}\left(C\left\|\vp_h^{k+1}\right\|_{H^1}^3 + \left\|\vp_h^k\right\|\right) \left\|\Delta_h\vp_h^{k+1}\right\| + \frac{1}{\gamma^2\epsilon^2}\left\|\mu_h^{k+1}\right\|^2 + \frac{1}{4}\left\|\Delta_h\vp_h^{k+1}\right\|^2 \non\\
&\leq& \frac{C}{\epsilon^4} + \frac{1}{\gamma^2\epsilon^2}\left\|\mu_h^{k+1}\right\|^2 + \frac{1}{2}\left\|\Delta_h\vp_h^{k+1}\right\|^2.
\end{eqnarray}
Therefore, we get 
\be
\left\|\Delta_h\vp_h^{k+1}\right\|^2 \leq \frac{2}{\gamma^2\epsilon^2}\left\|\mu_h^{k+1}\right\|^2 + \frac{2C}{\epsilon^4},\label{lemma_est_1_done}
\ee
which in turn proves \eqref{lemma_est_1}. Likewise by taking $\phi=\mu_h^{k+1}$ in \eqref{PD_CH2},  one derives
\begin{eqnarray}
\left\|\mu_h^{k+1}\right\|^2 &=& \frac{\gamma}{\epsilon}\left(f(\vp_h^{k+1},\vp_h^k), \mu_h^{k+1}\right) + \gamma\epsilon\left(\nabla\vp_h^{k+1},\nabla\mu_h^{k+1}\right) \non\\
&\leq& \frac{\gamma}{\epsilon}\left\|f(\vp_h^{k+1},\vp_h^k)\right\| \left\|\mu_h^{k+1}\right\| + \gamma\epsilon\left\|\nabla\vp_h^{k+1}\right\| \left\|\nabla\mu_h^{k+1}\right\| \non\\
&\leq& \frac{\gamma^2}{2\epsilon^2}\left\|f(\vp_h^{k+1},\vp_h^k)\right\|^2 + \frac{1}{2}\left\|\mu_h^{k+1}\right\|^2 + \frac{\gamma^2\epsilon^2}{2}\left\|\nabla\vp_h^{k+1}\right\|^2 + \frac{1}{2}\left\|\nabla\mu_h^{k+1}\right\|^2 \non\\
&\leq& \frac{\gamma^2}{2\epsilon^2}\left(\left\|\vp_h^{k+1}\right\|_{L^6}^3+\left\|\vp_h^k\right\|\right)^2 + \frac{1}{2}\left\|\mu_h^{k+1}\right\|^2 + \frac{\gamma^2\epsilon^2}{2}\left\|\nabla\vp_h^{k+1}\right\|^2 + \frac{1}{2}\left\|\nabla\mu_h^{k+1}\right\|^2 \non\\
&\leq& \frac{1}{2}\left\|\mu_h^{k+1}\right\|^2 + \frac{1}{2}\left\|\nabla\mu_h^{k+1}\right\|^2 + \frac{C\gamma^2}{2\epsilon^2} + \frac{C\gamma^2\epsilon^2}{2} . 
\end{eqnarray}
As a result, inequality  \eqref{lemma_est_2} holds, i.e.
\be
\left\|\mu_h^{k+1}\right\|^2 \leq \left\|\nabla\mu_h^{k+1}\right\|^2 + \frac{C\gamma^2}{\epsilon^2} + C\gamma^2\epsilon^2.\label{lemma_est_2_done}
\ee
Moreover, the inequality \eqref{lemma_est_3} follows from  \eqref{lemma_est_1}, \eqref{lemma_est_2} and \eqref{energy_law}. By Lemma \ref{nirenberg}, one has
\begin{eqnarray}
\left\|\varphi_h^{k+1}\right\|_{{\infty}} &\leq& C\left\|\Delta_h\varphi_h^{k+1}\right\|^{\frac{d}{2(6-d)}}\left\|\varphi_h^{k+1}\right\|_{L^6}^{\frac{3(4-d)}{2(6-d)}} + C\left\|\varphi_h^{k+1}\right\|_{L^6} \non\\
&\leq& C\left\|\Delta_h\varphi_h^{k+1}\right\|^{\frac{d}{2(6-d)}} + C.
\end{eqnarray}
Thus, an application of Young's inequality gives
\begin{eqnarray}
\left\|\varphi_h^{k+1}\right\|_{{\infty}}^{\frac{4(6-d)}{d}} \leq \left(C\left\|\Delta_h\varphi_h^{k+1}\right\|^{\frac{d}{2(6-d)}} + C\right)^{\frac{4(6-d)}{d}} 
 \leq \left(C\left\|\Delta_h\varphi_h^{k+1}\right\|^2 + C\right).\label{varphi_h^4_done}
\end{eqnarray}
Subsequently, a combination of \eqref{lemma_est_3}, \eqref{lemma_est_1_done}, and \eqref{varphi_h^4_done} yields \eqref{varphi_h^4}.

For the inequality \eqref{varphi_h^8}, we observe the following identity for any $v_h\in Y_h$, $\Delta_h v_h, \Delta_h^2 v_h \in \mathring{Y}_h$: 
\begin{eqnarray}
\left(\nabla v_h, \nabla\Delta_h^2 v_h\right) = \left\|\nabla\Delta_h v_h\right\|^2 = \left\|\Delta_h^2 v_h\right\|_{-1,h}^2, 
\end{eqnarray}
and that
\begin{eqnarray}
\left\|\left(\vp_h^{k+1}\right)^3 - \vp_h^k\right\|_{H^1}^2 &=& \left\|\left(\vp_h^{k+1}\right)^3 - \vp_h^k\right\|^2 + \left\|\nabla\left(\left(\vp_h^{k+1}\right)^3 - \vp_h^k\right)\right\|^2 \non\\
&\leq& 2\left\|\left(\vp_h^{k+1}\right)^3\right\|^2 + 2\left\|\vp_h^k\right\|^2 + 2\left\|\nabla\left(\vp_h^{k+1}\right)^3\right\|^2 + 2\left\|\nabla\vp_h^k\right\|^2 \non\\
&=& 2\left\|\left(\vp_h^{k+1}\right)\right\|_{L^6}^6 + 2\left\|\vp_h^k\right\|_{H^1}^2 + 2\left\|3\left(\vp_h^{k+1}\right)^2\nabla\vp_h^{k+1}\right\|^2 \non\\
&\leq& C\left\|\left(\vp_h^{k+1}\right)\right\|_{H^1}^6 + 2\left\|\vp_h^k\right\|_{H^1}^2 + 6\left\|\vp_h^{k+1}\right\|_{\infty}^4\left\|\nabla\vp_h^{k+1}\right\|^2 \non\\
&\leq&  C\left\|\left(\vp_h^{k+1}\right)\right\|_{H^1}^6 + 2\left\|\vp_h^k\right\|_{H^1}^2 + 6\left\|\nabla\vp_h^{k+1}\right\|^2 \left(\frac{d}{6-d}\left\|\vp_h^{k+1}\right\|_{\infty}^{\frac{4(6-d)}{d}} + \frac{6-2d}{6-d}\right) \non\\
&\leq& C\left\|\vp_h^{k+1}\right\|_{\infty}^{\frac{4(6-d)}{d}} + C.
\end{eqnarray}
Then by taking $\phi_h=\Delta_h^2 \vp_h^{k+1}$ in \eqref{PD_CH2}, one obtains 
\begin{eqnarray}
\left\|\nabla\Delta_h\vp_h^{k+1}\right\|^2 &=& \frac{1}{\gamma\epsilon}\left(\mu_h^{k+1}, \Delta_h^2\vp_h^{k+1}\right) - \frac{1}{\epsilon^2}\left(\left(\vp_h^{k+1}\right)^3 - \vp_h^k, \Delta_h^2\vp_h^{k+1}\right) \non\\
&\leq& -\frac{1}{\gamma\epsilon}\left(\nabla\mu_h^{k+1}, \nabla\Delta_h\vp_h^{k+1}\right) + \frac{1}{\epsilon^2}\left\|\left(\vp_h^{k+1}\right)^3 - \vp_h^k \right\|_{H^1} \left\|\Delta_h^2\vp_h^{k+1}\right\|_{-1,h} \non\\
&\leq& \frac{1}{\gamma\epsilon}\left\|\nabla\mu_h^{k+1}\right\| \left\|\nabla\Delta_h\vp_h^{k+1}\right\| + \frac{1}{\epsilon^2}\left\|\left(\vp_h^{k+1}\right)^3 - \vp_h^k \right\|_{H^1} \left\|\nabla\Delta_h\vp_h^{k+1}\right\| \non\\
&\leq& \frac{1}{\gamma^2\epsilon^2}\left\|\nabla\mu_h^{k+1}\right\|^2 + \frac{1}{\epsilon^4}\left\|\left(\vp_h^{k+1}\right)^3 - \vp_h^k \right\|_{H^1}^2 + \frac{1}{2}\left\|\nabla\Delta_h\vp_h^{k+1}\right\|^2 \non\\
&\leq& \frac{1}{\gamma^2\epsilon^2}\left\|\nabla\mu_h^{k+1}\right\|^2 + \frac{C}{\epsilon^4}\left\|\vp_h^{k+1}\right\|_{\infty}^{\frac{4(6-d)}{d}} + \frac{C}{\epsilon^4} + \frac{1}{2}\left\|\nabla\Delta_h\vp_h^{k+1}\right\|^2,
\end{eqnarray}
which yields that
\begin{eqnarray}
\left\|\nabla\Delta_h\vp_h^{k+1}\right\|^2 \leq \frac{2}{\gamma^2\epsilon^2}\left\|\nabla\mu_h^{k+1}\right\|^2 + \frac{C}{\epsilon^4}\left\|\vp_h^{k+1}\right\|_{\infty}^{\frac{4(6-d)}{d}} + \frac{C}{\epsilon^4}. \label{nabla_delta_done}
\end{eqnarray}
Also notice that $\left(\vp_h^{k},1\right)\equiv \left(\vp_h^0,1\right) = C,\ \forall 0\leq k \leq K$, by taking $v_h=1$ in \eqref{PD_CH1}. By Lemma \ref{nirenberg}, we derive
\begin{eqnarray}
\left\|\vp_h^{k+1}\right\|_{\infty} &\leq& \left\|\vp_h^{k+1}-\overline{\vp_h^{k+1}}\right\|_{\infty} + \left|\overline{\vp_h^{k+1}}\right| \leq C \left\|\nabla\Delta_h\vp_h^{k+1}\right\|^{\frac{d}{4(6-d)}} \left\|\nabla\vp_h^{k+1}\right\|^{\frac{24-5d}{4(6-d)}} 
  + C\left\|\nabla\vp_h^{k+1}\right\| + \left|\overline{\vp_h^{0}}\right| \non\\
&\leq& C\left\|\nabla\Delta_h \vp_h^{k+1}\right\|^{\frac{d}{4(6-d)}} + C,
\end{eqnarray}
so that 
\begin{eqnarray}
\left\|\vp_h^{k+1}\right\|_{\infty}^{\frac{8(6-d)}{d}} \leq C\left\|\nabla\Delta_h \vp_h^{k+1}\right\|^2 + C. \label{varphi_h^8_done}
\end{eqnarray}
Combining \eqref{nabla_delta_done}, \eqref{varphi_h^8_done}, \eqref{energy_law} and \eqref{varphi_h^4}, one readily derives \eqref{varphi_h^8}. This completes the proof. \hfill 
\end{proof}

\section{The optimal rate error analysis}\label{Err-ana}
In this section we provide a convergence analysis and error estimate for the numerical scheme \eqref{PD_CH1}--\eqref{initialC}.  Further regularity assumptions for the weak solution are needed in the analysis. 

\begin{myassum}\label{higher_regularities}
We assume that  weak solutions to the CHSD system \eqref{weak_CH1}-\eqref{weak_D} have the following additional regularities
\begin{eqnarray}
&& \varphi \in L^{\infty}\left(0,T;W^{1,6}(\Omega)\right)\bigcap L^4\left(0,T;H^1(\Omega)\right)\bigcap H^{2}\left(0,T;L^2(\Omega)\right) \bigcap L^{\infty}\left(0,T;H^{q+1}(\Omega)\right), \\
&& \mu \in L^{\infty}\left(0,T;H^{q+1}(\Omega)\right), \\
&& \ub_c \in L^{\infty}\left(0,T;\left[H^{q+1}(\Omega_c)\right]^d\right)\bigcap  W^{1,4}\left(0,T;\left[L^{2}(\Omega_c)\right]^d\right) \bigcap  H^{2}\left(0,T;\left[L^{2}(\Omega_c)\right]^d\right), \\
&& \ub_m \in L^{\infty}\left(0,T;\left[H^{q+1}(\Omega_m)\right]^d\right)\bigcap W^{1,4}\left(0,T;\left[L^2(\Omega_m)\right]^d\right)\bigcap H^{2}\left(0,T;\left[L^2(\Omega_m)\right]^d\right),
\end{eqnarray}
where $q\geq1$ is the spatial approximation order.
\end{myassum}

The following assumptions are also made, on the the parameters of the problem
\begin{eqnarray}\label{para_assump}
&M_0\leq{\rm M(\vp)}\leq M_1, \quad  |M^\prime| \leq C, \quad 
\nu_0\leq\nu(\vp)\leq\nu_1, \quad |\nu^\prime| \leq C.
\end{eqnarray}


For  the weak solution $(\ub_c, P_c, \mathbf{u}_m, P_m, \vp, \mu)$ to the CHSD system \eqref{weak_CH1}--\eqref{weak_D}, we set
\begin{eqnarray}
&& \rho^{\vp}(x,t):=\vp(x,t)-\mathcal{P}\vp(x,t), \ \ \ \ \ \ \ \rho^{\mu}(x,t):=\mu(x,t)-\widetilde{\mathcal{P}}^{\vp(t)}\mu(x,t), \label{rho_1}\\
&& \rho^{\ub}(x,t)\Big|_{\Omega_j} = \rho_{j}^{\ub}(x,t):=\ub_j(x,t)-\mathcal{P}_{j,u}^{\vp(t)}\ub_j(x,t),\ \ \ \ j\in\{c,m\}, \label{rho_2}
\end{eqnarray}
specially, for $0 \leq k \leq K,\ j\in\{c,m\}$,
\begin{eqnarray}
&& \rho^{\varphi,k}\Big|_{\Omega_j} = \rho_j^{\varphi,k}:=\left(\varphi^k-\mathcal{P}\varphi^k\right)\Big|_{\Omega_j},\ \ \ \rho^{\mu,k}\Big|_{\Omega_j}= \rho_j^{\mu,k}:=\left(\mu^k-\widetilde{\mathcal{P}}^k\mu^k\right)\Big|_{\Omega_j}, \label{rho_3}\\
&& \rho^{\mathbf{u},k}\Big|_{\Omega_j} = \rho_j^{\mathbf{u},k}:=\mathbf{u}_j^k-\mathcal{P}_{j,u}^k\mathbf{u}_j^k,\ \ \ \rho^{p,k}\Big|_{\Omega_j} = \rho_j^{p,k}:=P_j^k-\mathcal{P}_{j,p}^k P_j^k, \label{rho_4}
\end{eqnarray}
and for $0 \leq k \leq K-1,\ j\in\{c,m\}$,
\begin{eqnarray}
&& R^{\varphi,k+1}\Big|_{\Omega_j} = R_j^{\varphi,k+1} := \left(\delta_t \mathcal{P}\varphi^{k+1}-\partial_t\varphi^{k+1}\right)\Big|_{\Omega_j},\ \ \ \ R^{\ub,k+1}\Big|_{\Omega_j} = R_j^{\mathbf{u},k+1}:=\delta_t\mathcal{P}_{j,u}^{k+1} \mathbf{u}_j^{k+1}-\partial_t\mathbf{u}_j^{k+1}, \non\\
&& R^{k+1} := \left\|\vp^{k+1}-\vp^k\right\|_{H^1}^2 + \left\|\ub_c^{k+1}-\ub_c^{k}\right\|^2 + \left\|\ub_m^{k+1}-\ub_m^{k}\right\|^2 = \left\|\vp^{k+1}-\vp^k\right\|_{H^1}^2 + \left\|\ub^{k+1}-\ub^{k}\right\|^2.
\end{eqnarray}

The error functions are defined as follows, for $j\in\{c,m\}$ and $0\leq k\leq K$:
\begin{eqnarray}
&& \sigma^{\vp,k}\Big|_{\Omega_j} = \sigma_j^{\vp,k}:=\left(\mathcal{P}\vp^k-\vp_h^k\right)\Big|_{\Omega_j},\ \ \ \ e^{\vp,k}\Big|_{\Omega_j} = e_j^{\vp,k}:=\left(\vp^k-\vp_h^k\right)\Big|_{\Omega_j},\label{ef_1}\\
&& \sigma^{\mu,k}\Big|_{\Omega_j} = \sigma_j^{\mu,k}:=\left(\widetilde{\mathcal{P}}^k\mu^k-\mu_h^k\right)\Big|_{\Omega_j},\ \ \ \ e^{\mu,k}\Big|_{\Omega_j} = e_j^{\mu,k} := \left(\mu^k - \mu_h^k\right)\Big|_{\Omega_j}, \\
&& \sigma^{\ub,k}\Big|_{\Omega_j} = \sigma_j^{\ub,k}:=\mathcal{P}_{j,u}^k \ub_j^k-\ub_{j,h}^k,\ \ \ \  e^{\ub,k}\Big|_{\Omega_j} = e_j^{\ub,k}:=\ub_j^k-\ub_{j,h}^k,\\
&& \sigma^{p,k}\Big|_{\Omega_j} = \sigma_j^{p,k}:=\mathcal{P}_{j,p}^k P_j^k-P_{j,h}^k,\ \ \ \  e^{p,k}\Big|_{\Omega_j} = e_j^{p,k}:=P_j^k-P_{j,h}^k. \label{ef_4}
\end{eqnarray}
Note that the numerical solution $\vp_h^k$ satisfies mass-conservation by choosing $v_h=1$ in \eqref{PD_CH1}, same as the weak solution $\vp$. Recall also that $\vp_h^0=\mathcal{P}\vp^0$. Then by the definition of Ritz projection we see that $\left(\vp^k,1\right) = \left(\mathcal{P}\vp^k,1\right) = \left(\vp_h^k,1\right) \equiv C_0$ for all $0\leq k\leq K$. This enables one to apply Poincar\'{e} inequality to $\rho^{\vp,k}, \sigma^{\vp,k}, e^{\vp,k}, \delta_t \sigma^{\vp,k+1}$ for $0\leq k \leq K$. We shall also make use of the fact that $\sigma^{\varphi,k},\ \delta_t \sigma^{\varphi,k+1} \in \mathring{Y}_h$.

Given any $t\in[0,T]$, the solution to the CHSD system satisfies
\begin{subequations}
\begin{eqnarray}
&&(\delta_t \mathcal{P}\varphi^{k+1},v) + ({\rm M}(\vp^{k+1})\nabla \widetilde{\mathcal{P}}^{k+1} \mu^{k+1},\nabla v) - (\ub^{k+1} \varphi^{k+1} ,\nabla v) = ( R^{\varphi,k+1},v), \label{ww_CH1}\\
&&\gamma\epsilon(\nabla \mathcal{P}\varphi^{k+1},\nabla \phi) -(\widetilde{\mathcal{P}}^{k+1}\mu^{k+1},\phi) + \frac{\gamma}{\epsilon}(f(\varphi^{k+1}),\phi) = (\rho^{\mu,k+1},\phi), \label{ww_CH2}\\
&&\rho_0(\delta_t\mathcal{P}_{c,u}^{k+1}\ub_c^{k+1},\mathbf{v}_c )_c + a_c(\mathcal{P}_{c,u}^{k+1}\ub_c^{k+1},\mathbf{v}_c) + b_c(\mathbf{v}_c, \mathcal{P}_{c,p}^{k+1} P_c^{k+1}) + \int_{\Gamma_{cm}} \mathcal{P}_{m,p}^{k+1}P_m^{k+1} (\mathbf{v}_c\cdot \mathbf{n}_{cm})dS \non\\
&&\ \ \ \ -b_c (\mathcal{P}_{c,u}^{k+1}\ub_c^{k+1}, q_c) +(\vp_c^{k+1} \nabla \mu_c^{k+1},  \mathbf{v}_c )_c = \rho_0( R_c^{\ub,k+1},\mathbf{v}_c)_c, \label{ww_S}\\
&&\frac{\rho_0}{\chi}(\delta_t\mathcal{P}_{m,u}^{k+1}\ub_m^{k+1},\mathbf{v}_m )_m + a_m(\mathcal{P}_{m,u}^{k+1}\ub_m^{k+1},  \mathbf{v}_m)_m + b_m(\mathbf{v}_m,\mathcal{P}_{m,p}^{k+1}P_m^{k+1}) - b_m(\mathcal{P}_{m,u}^{k+1}\ub_m^{k+1}, q_m)\non \\
&&\ \ \ \ -\int_{\Gamma_{cm}} \mathcal{P}_{c,u}^{k+1}\ub_c^{k+1} \cdot \mathbf{n}_{cm} q_m\, dS +(\vp_m^{k+1} \nabla \mu_m^{k+1},  \mathbf{v}_m )_m= \frac{\rho_0}{\chi}( R_m^{\ub,k+1},\mathbf{v}_m )_m, \label{ww_D}
\end{eqnarray}
\end{subequations}
for all $v,\ \phi\in Y_h,\ \mathbf{v}_j\in  \mathbf{X}_j^h,\ q_j\in  M_j^h,\ j\in\{c,m\}$ and  $0\leq k\leq K-1$.

Subtracting \eqref{PD_CH1}-\eqref{PD_D} from \eqref{ww_CH1}-\eqref{ww_D}, we obtain 
\begin{subequations}
\begin{eqnarray}
&&(\delta_t\sigma^{\varphi,k+1},v) + ({\rm M}(\vp_h^k)\nabla \sigma^{\mu,k+1},\nabla v) = -\left(({\rm M}(\vp^{k+1}) - {\rm M}(\vp_h^k))\nabla \widetilde{\mathcal{P}}^{k+1}\mu^{k+1}, \nabla v\right) \nonumber\\
&&\ \ \ \ \ \ + (\ub^{k+1} \varphi^{k+1} - \baru^{k+1}_h \varphi^k_h ,\nabla v) + ( R^{\varphi,k+1},v), \label{pd_err_CH1}\\
&&\gamma\epsilon(\nabla \sigma^{\varphi,k+1},\nabla \phi) -(\sigma^{\mu,k+1},\phi) = (\rho^{\mu,k+1},\phi) - \frac{\gamma}{\epsilon}\left(f(\varphi^{k+1})-f(\varphi^{k+1}_h, \varphi^k_h),\phi\right), \label{pd_err_CH2}\\
&&\rho_0(\delta_t\sigma_c^{\ub,k+1},\mathbf{v}_c )_c + \int_{\Gamma_{cm}} \sigma_m^{p,k+1} (\mathbf{v}_c\cdot \mathbf{n}_{cm})\,ds\non\\
&&\ \ \ \ \ \ + a_c^k(\sigma_c^{\ub,k+1},\mathbf{v}_c) + b_c(\mathbf{v}_c, \sigma_c^{p,k+1}) -b_c (\sigma_c^{\ub,k+1}, q_c) \non\\
&&\ \ \ \ \ = \rho_0( R_c^{\ub,k+1},\mathbf{v}_c)_c - \left(\vp_c^{k+1}\nabla\mu_c^{k+1}-\vp_{c,h}^k \nabla \mu_{c,h}^{k+1},  \mathbf{v}_c \right)_c \non\\
&&\ \ \ \ \ \ - 2\left((\nu(\vp_c^{k+1})-\nu(\vp_{c,h}^k))\mathbb{D} (\mathcal{P}_{c,u}^{k+1}\mathbf{u}_c^{k+1}),\mathbb{D}(\mathbf{v}_{c})\right)_c \non\\
&&\ \ \ \ \ \ + \sum_{i=1}^{d-1}\int_{\Gamma_{cm}}\alpha_{BJSJ}\frac{\nu(\vp_c^{k+1})-\nu(\vp_{c,h}^k)}{\sqrt{{\rm tr} (\Pi)}}
(\mathcal{P}_{c,u}^{k+1}\mathbf{u}_c^{k+1}\cdot\btau_i)(\mathbf{v}_{c}\cdot\btau_i)dS , \label{pd_err_S}\\
&&\frac{\rho_0}{\chi}(\delta_t\sigma_m^{\ub,k+1},\mathbf{v}_m )_m - \int_{\Gamma_{cm}} \sigma_c^{\ub,k+1} \cdot \mathbf{n}_{cm} q_m\ ds + a_m^k(\sigma_m^{\ub,k+1},  \mathbf{v}_m)_m \non\\
&&\ \ \ \ \ \ + b_m(\mathbf{v}_m,\sigma_m^{p,k+1}) - b_m(\sigma_m^{\ub,k+1}, q_m)\non \\
&&\ \ \ \ \ = \frac{\rho_0}{\chi}( R_m^{\ub,k+1},\mathbf{v}_m )_m - \left(\vp_m^{k+1}\nabla\mu_m^{k+1}-\vp_{m,h}^k \nabla \mu_{m,h}^{k+1},  \mathbf{v}_m \right)_m \non\\
&&\ \ \ \ \ \ - \left((\nu(\vp_m^{k+1})-\nu(\varphi_{m,h}^{k}))\Pi^{-1}\mathcal{P}_{m,u}^{k+1} \mathbf{u}_{m}^{k+1}, \mathbf{v}_{m,h}\right)_m ,  \label{pd_err_D}
\end{eqnarray}
\end{subequations}
for all $0\leq k\leq K-1$, $v,\ \phi\in Y_h,\ \mathbf{v}_j\in  \mathbf{X}_j^h,\ q_j\in  M_j^h,\ j\in\{c,m\}$.

Setting $v=\sigma^{\mu,k+1}$ in \eqref{pd_err_CH1},\ \ $\phi=\delta_t\sigma^{\phi,k+1}$ in \eqref{pd_err_CH2},\ \ $\vb_c=\sigma_c^{\ub,k+1},\ q_c=\sigma_c^{p,k+1}$ in \eqref{pd_err_S},\ \ $\vb_m=\sigma_m^{\ub,k+1},\ q_m=\sigma_m^{p,k+1}$ in \eqref{pd_err_D}, adding the resulting equations, and noticing that for $d=2,3$,
\begin{eqnarray}
&& M_0\leq{\rm M(\vp)}\leq M_1,\ \ \ \ \ \ \nu_0\leq\nu(\vp)\leq\nu_1,\ \ \ \ \ \ \lambda_{\max}(\Pi) \leq \lambda,\ \ \ \ \ \ {\rm tr}(\Pi) \leq d\lambda, \non\\
&& \left\|\ub\right\|^2 = \left\|\Pi^{1/2}\Pi^{-1/2}\ub\right\|^2 \leq \left\|\Pi^{1/2}\right\|_2^2 \left\|\Pi^{-1/2}\ub\right\|^2 = \lambda_{\max}(\Pi)\left\|\Pi^{-1/2}\ub\right\|^2 \leq \lambda \left\|\Pi^{-1/2}\ub\right\|^2,
\end{eqnarray}
we derive the following error equation for the numerical scheme:
\begin{eqnarray}\label{pd_err}
&& M_0\left\|\nabla\sigma^{\mu,k+1}\right\|^2 + \frac{\gamma\epsilon}{2\tau}\left(\left\|\nabla\sigma^{\varphi,k+1}\right\|^2 - \left\|\nabla\sigma^{\varphi,k}\right\|^2 + \left\|\nabla(\sigma^{\varphi,k+1}-\sigma^{\varphi,k})\right\|^2\right)\nonumber\\
&&\ \ +\ \frac{\rho_0}{2\tau}\left(\left\|\sigma_c^{\ub,k+1}\right\|^2 - \left\|\sigma_c^{\ub,k}\right\|^2 + \left\|\sigma_c^{\ub,k+1}-\sigma_c^{\ub,k}\right\|^2\right) + \alpha_{BJSJ}\frac{\nu_0}{\sqrt{d\lambda}}\sum_{i=1}^{d-1}\left\|\sigma_c^{\ub,k+1}\cdot\btau_i\right\|_{cm}^2 \non\\
&&\ \ +\ 2\nu_0\left\|\mathbb{D}(\sigma_c^{\ub,k+1})\right\|^2 + \frac{\rho_0}{2\tau\chi}\left(\left\|\sigma_m^{\ub,k+1}\right\|^2 - \left\|\sigma_m^{\ub,k}\right\|^2 + \left\|\sigma_m^{\ub,k+1}-\sigma_m^{\ub,k}\right\|^2\right) + \frac{\nu_0}{\lambda}\left\|\sigma_m^{\ub,k+1}\right\|^2 \nonumber\\
&&=\ - \left(\left({\rm M}(\vp^{k+1})-{\rm M}(\vp_h^k)\right)\nabla\widetilde{\mathcal{P}}^{k+1}\mu^{k+1}, \nabla\sigma^{\mu,k+1}\right) \non\\
&&\ \ - 2\left(\left(\nu(\vp_c^{k+1})-\nu(\vp_{c,h}^k)\right)\mathbb{D}(\mathcal{P}_{c,u}^{k+1} \ub_c^{k+1}), \mathbb{D}(\sigma_c^{\ub,k+1})\right)_c \non\\
&&\ \ - \sum_{i=1}^{d-1}\int_{\Gamma_{cm}}\alpha_{BJSJ}\frac{\nu(\vp_c^{k+1})-\nu(\vp_{c,h}^k)}{\sqrt{{\rm tr}(\Pi)}} \left(\mathcal{P}_{c,u}^{k+1} \ub_c^{k+1}\cdot\btau_i\right) \left(\sigma_c^{\ub,k+1}\cdot\btau_i\right)dS \non\\
&&\ \ - \left(\left(\nu(\vp_m^{k+1})-\nu(\vp_{m,h}^k)\right)\Pi^{-1}\mathcal{P}_{m,u}^{k+1} \ub_{m}^{k+1}, \sigma_m^{\ub,k+1}\right)_m \non\\
&&\ \ + \frac{\rho_0}{\chi}\left( R_m^{\ub,k+1},\sigma_m^{\ub,k+1}\right)_m + \rho_0\left( R_c^{\ub,k+1},\sigma_c^{\ub,k+1}\right)_c + \left( R^{\varphi,k+1},\sigma^{\mu,k+1}\right) \non\\
&&\ \  + \left(\rho^{\mu,k+1},\delta_t\sigma^{\varphi,k+1}\right) + \left(\ub^{k+1}\varphi^{k+1}-\overline{\ub}_h^{k+1}\varphi_h^k,\nabla\sigma^{\mu,k+1}\right)\nonumber\\
&&\ \ - \frac{\gamma}{\epsilon}\left(f(\varphi^{k+1},\varphi^{k+1})-f(\varphi_h^{k+1},\varphi_h^k),
\delta_t\sigma^{\varphi,k+1}\right) - \left(\varphi^{k+1}\nabla\mu^{k+1}-\varphi_{h}^k\nabla\mu_{h}^{k+1},\sigma^{\ub,k+1}\right)\\
&&:=\sum_{j=1}^{11} I_j \nonumber 
\end{eqnarray}
where we have designated the eleven terms on the right-hand side of  \eqref{pd_err} by $I_j, j=1,2 \cdots 11$.
Now we estimate the $I_j$s in a series of lemmas.
\begin{lemma}[Estimate of the first term $I_1$]\label{lemma_RHS1}
Suppose $(\varphi, \mu, \ub_c, \ub_m, P_c, P_m)$ is a weak solution to \eqref{ww_CH1}--\eqref{ww_D} with the additional regularities described in Assumption \ref{higher_regularities}, $d=2,3$. Set $M_0$ as the lower bound of the mobility ${\rm M}(\vp)$. Then the first term $I_1$ of  RHS  of \eqref{pd_err} satisfies
\begin{eqnarray}\label{RHS1}
\left|- \left(\left({\rm M}(\vp^{k+1})-{\rm M}(\vp_h^k)\right)\nabla\widetilde{\mathcal{P}}^{k+1}\mu^{k+1}, \nabla\sigma^{\mu,k+1}\right)\right| \leq C\left(R^{k+1} + \left\|\nabla e^{\vp,k}\right\|^2\right) + \frac{M_0}{12}\left\|\nabla\sigma^{\mu,k+1}\right\|^2,
\end{eqnarray}
for a constant $C$  independent of $\tau$ and $h$.
\end{lemma}
\begin{proof}
We split the term into two parts as follows
\begin{eqnarray}
- \left(\left({\rm M}(\vp^{k+1})-{\rm M}(\vp_h^k)\right)\nabla\widetilde{\mathcal{P}}^{k+1}\mu^{k+1}, \nabla\sigma^{\mu,k+1}\right) &=& \left(\left({\rm M}(\vp^{k+1})-{\rm M}(\vp_h^k)\right)\nabla\rho^{\mu,k+1}, \nabla\sigma^{\mu,k+1}\right) \non\\
&& - \left(\left({\rm M}(\vp^{k+1})-{\rm M}(\vp_h^k)\right)\nabla\mu^{k+1}, \nabla\sigma^{\mu,k+1}\right).
\end{eqnarray}
By the inverse inequality, there exists a constant $\theta_1>0$ such that for all $0\leq k\leq K-1$, we have
\begin{eqnarray}
&& \left|\left(\left({\rm M}(\vp^{k+1})-{\rm M}(\vp_h^k)\right)\nabla\rho^{\mu,k+1}, \nabla\sigma^{\mu,k+1}\right)\right| \non\\
&\leq& C\left\|{\rm M}(\vp^{k+1})-{\rm M}(\vp_h^k)\right\|_6 \left\|\nabla\rho^{\mu,k+1}\right\| \left\|\nabla\sigma^{\mu,k+1}\right\|_3 \non\\
&\leq& C\left\|\vp^{k+1}-\vp_h^k\right\|_6 h\left\|\mu^{k+1}\right\|_{H^2} h^{d/3-d/2}\left\|\nabla\sigma^{\mu,k+1}\right\| \non\\
&\leq& Ch^{1-d/6}\left\|\vp^{k+1}-\vp_h^k\right\|_{H^1} \left\|\nabla\sigma^{\mu,k+1}\right\| \non\\
&\leq& \frac{C}{\theta_1}\left\|\vp^{k+1}-\vp_h^k\right\|_{H^1}^2 + \frac{\theta_1}{2}\left\|\nabla\sigma^{\mu,k+1}\right\|^2 \non\\
&\leq& \frac{C}{\theta_1}\left(\left\|\vp^{k+1}-\vp^k\right\|_{H^1}^2 + \left\|e^{\vp,k}\right\|_{H^1}^2\right) + \frac{\theta_1}{2}\left\|\nabla\sigma^{\mu,k+1}\right\|^2 \non\\
&\leq& \frac{C}{\theta_1}\left(R^{k+1} + \left\|\nabla e^{\vp,k}\right\|^2\right) + \frac{\theta_1}{2}\left\|\nabla\sigma^{\mu,k+1}\right\|^2, \label{split-1}
\end{eqnarray}
and similarly,
\begin{eqnarray}
&& \left|\left(\left({\rm M}(\vp^{k+1})-{\rm M}(\vp_h^k)\right)\nabla\mu^{k+1}, \nabla\sigma^{\mu,k+1}\right)\right| \non\\
&\leq& C\left\|{\rm M}(\vp^{k+1})-{\rm M}(\vp_h^k)\right\|_6 \left\|\nabla\mu^{k+1}\right\|_3 \left\|\nabla\sigma^{\mu,k+1}\right\| \non\\
&\leq& C\left\|\vp^{k+1}-\vp_h^k\right\|_6 \left\|\nabla\sigma^{\mu,k+1}\right\| \non\\
&\leq& C\left\|\vp^{k+1}-\vp_h^k\right\|_{H^1} \left\|\nabla\sigma^{\mu,k+1}\right\| \non\\
&\leq& \frac{C}{\theta_1}\left(R^{k+1} + \left\|\nabla e^{\vp,k}\right\|^2\right) + \frac{\theta_1}{2}\left\|\nabla\sigma^{\mu,k+1}\right\|^2. \label{split-2}
\end{eqnarray}
Combining \eqref{split-1} and \eqref{split-2} and  choosing $\theta_1 = \frac{M_0}{12}$, one obtains \eqref{RHS1}. This completes the proof.\hfill 
\end{proof}

The estimates of $I_2, I_3, I_4$  in \eqref{pd_err} are summarized in the following lemma.
\begin{lemma}[Estimates of $I_2, I_3, I_4$ ]\label{lemma_RHS2-4}
The assumptions are the same as in Lemma \ref{lemma_RHS1}.  Then  $I_2, I_3, I_4$ of RHS in \eqref{pd_err} satisfy
\begin{eqnarray}
&& \left|- 2\left(\left(\nu(\vp_c^{k+1})-\nu(\vp_{c,h}^k)\right)\mathbb{D}(\mathcal{P}_{c,u}^{k+1} \ub_c^{k+1}), \mathbb{D}(\sigma_c^{\ub,k+1})\right)_c\right| \non\\
&&\ \ \ \ \ \ \leq C\left(R^{k+1} + \left\|\nabla e^{\vp,k}\right\|^2 \right) + \frac{\nu_0}{2}\left\|\mathbb{D}(\sigma_c^{\ub,k+1})\right\|^2, \label{RHS2}\\
&& \left|- \sum_{i=1}^{d-1}\int_{\Gamma_{cm}}\alpha_{BJSJ}\frac{\nu(\vp_c^{k+1})-\nu(\vp_{c,h}^k)}{\sqrt{{\rm tr}(\Pi)}} \left(\mathcal{P}_{c,u}^{k+1} \ub_c^{k+1}\cdot\btau_i\right) \left(\sigma_c^{\ub,k+1}\cdot\btau_i\right)dS\right| \non\\
&&\ \ \ \ \ \ \leq C\left( R^{k+1} + \left\|\nabla e^{\vp,k}\right\|^2 \right) + \alpha_{BJSJ}\frac{\nu_0}{2\sqrt{d\lambda}}\sum_{i=1}^{d-1}\left\|\sigma_c^{\ub,k+1}\cdot\btau_i\right\|_{cm}^2, \label{RHS3}\\
&& \left|- \left(\left(\nu(\vp_m^{k+1})-\nu(\vp_{m,h}^k)\right)\Pi^{-1}\mathcal{P}_{m,u}^{k+1} \ub_{m,h}^{k+1}, \sigma_m^{\ub,k+1}\right)_m\right| \non\\
&&\ \ \ \ \ \ \leq C\left( R^{k+1} + \left\|\nabla e^{\vp,k}\right\|^2 \right) + \frac{\nu_0}{4\lambda}\left\|\sigma_m^{\ub,k+1}\right\|^2 \label{RHS4}, 
\end{eqnarray}
where $C$s are constants independent of $\tau$ and $h$.
\end{lemma}
\begin{proof}
The inequality \eqref{RHS2} is derived the same way as \eqref{RHS1}, that is
\begin{eqnarray}
&& \left|- 2\left(\left(\nu(\vp_c^{k+1})-\nu(\vp_{c,h}^k)\right)\mathbb{D}(\mathcal{P}_{c,u}^{k+1} \ub_c^{k+1}), \mathbb{D}(\sigma_c^{\ub,k+1})\right)_c\right| \non\\
&\leq& \left|2\left(\left(\nu(\vp_c^{k+1})-\nu(\vp_{c,h}^k)\right)\mathbb{D}(\rho_{c}^{\ub,k+1}), \mathbb{D}(\sigma_c^{\ub,k+1})\right)_c\right| + \left|2\left(\left(\nu(\vp_c^{k+1})-\nu(\vp_{c,h}^k)\right)\mathbb{D}(\ub_{c}^{k+1}), \mathbb{D}(\sigma_c^{\ub,k+1})\right)_c\right| \non\\
&\leq& 2\left\|\nu(\vp_c^{k+1})-\nu(\vp_{c,h}^k)\right\|_6 \left\|\mathbb{D}(\rho_c^{\ub,k+1})\right\| \left\|\mathbb{D}(\sigma_c^{\ub,k+1})\right\|_3 + 2\left\|\nu(\vp_c^{k+1})-\nu(\vp_{c,h}^k)\right\|_6 \left\|\mathbb{D}(\ub_c^{k+1})\right\|_3 \left\|\mathbb{D}(\sigma_c^{\ub,k+1})\right\| \non\\
&\leq& C h^{1-d/6}\left\|\vp_c^{k+1}-\vp_{c,h}^k\right\|_6 \left\|\ub_c^{k+1}\right\|_{H^2} \left\|\mathbb{D}(\sigma_c^{\ub,k+1})\right\| + C\left\|\vp_c^{k+1}-\vp_{c,h}^k\right\|_6 \left\|\mathbb{D}(\sigma_c^{\ub,k+1})\right\| \non\\
&\leq& C\left\|\vp_c^{k+1}-\vp_{c,h}^k\right\|_{H^1} \left\|\mathbb{D}(\sigma_c^{\ub,k+1})\right\| \non\\
&\leq& \frac{C}{\theta_2}\left(R^{k+1} + \left\|\nabla e^{\vp,k}\right\|^2\right) + \theta_2\left\|\mathbb{D}(\sigma_c^{\ub,k+1})\right\|^2.
\end{eqnarray}
With an application of Lemma \ref{trace}, one has
\begin{eqnarray}
&& \left|- \sum_{i=1}^{d-1}\int_{\Gamma_{cm}}\alpha_{BJSJ}\frac{\nu(\vp_c^{k+1})-\nu(\vp_{c,h}^k)}{\sqrt{{\rm tr}(\Pi)}} \left(\mathcal{P}_{c,u}^{k+1} \ub_c^{k+1}\cdot\btau_i\right) \left(\sigma_c^{\ub,k+1}\cdot\btau_i\right)dS\right| \non\\
&\leq& \sum_{i=1}^{d-1}C\left\|\nu(\vp_c^{k+1})-\nu(\vp_{c,h}^k)\right\|_{L^4(\Gamma_{cm})} \left(\left\|\rho_c^{\ub,k+1}\right\|_{L^4(\Gamma_{cm})} + \left\|\ub_c^{k+1}\right\|_{L^4(\Gamma_{cm})}\right) \left\|\sigma_c^{\ub,k+1}\cdot\btau_i\right\|_{cm} \non\\
&\leq& \sum_{i=1}^{d-1}C\left\|\vp_c^{k+1}-\vp_{c,h}^k\right\|_{L^4(\Gamma_{cm})} \left(\left\|\mathbb{D}(\rho_c^{\ub,k+1})\right\|_{L^2(\Omega_c)} + \left\|\mathbb{D}(\ub_c^{k+1})\right\|_{L^2(\Omega_c)}\right) \left\|\sigma_c^{\ub,k+1}\cdot\btau_i\right\|_{cm} \non\\
&\leq& \sum_{i=1}^{d-1}C\left\|\vp_c^{k+1}-\vp_{c,h}^k\right\|_{H^1(\Omega_c)} \left(h\left\|\ub_c^{k+1}\right\|_{H^2(\Omega_c)} + \left\|\ub_c^{k+1}\right\|_{H^2(\Omega_c)}\right) \left\|\sigma_c^{\ub,k+1}\cdot\btau_i\right\|_{cm} \non\\
&\leq& \sum_{i=1}^{d-1}C\left\|\vp^{k+1}-\vp_{h}^k\right\|_{H^1} \left\|\sigma_c^{\ub,k+1}\cdot\btau_i\right\|_{cm} \non\\
&\leq& \frac{C}{\theta_3}\left(R^{k+1} + \left\|\nabla e^{\vp,k}\right\|^2\right) + \theta_3\sum_{i=1}^{d-1}\left\|\sigma_c^{\ub,k+1}\cdot\btau_i\right\|_{cm}^2.
\end{eqnarray}
Likewise,
\begin{eqnarray}
&& \left|- \left(\left(\nu(\vp_m^{k+1})-\nu(\vp_{m,h}^k)\right)\Pi^{-1}\mathcal{P}_{m,u}^{k+1} \ub_{m,h}^{k+1}, \sigma_m^{\ub,k+1}\right)_m\right| \non\\
&\leq& C\left\|\nu(\vp_m^{k+1})-\nu(\vp_{m,h}^k)\right\|_6 \left(\left\|\rho_m^{\ub,k+1}\right\| \left\|\sigma_m^{\ub,k+1}\right\|_3 + \left\|\ub_m^{k+1}\right\|_3 \left\|\sigma_m^{\ub,k+1}\right\|\right) \non\\
&\leq& C\left\|\vp_m^{k+1}-\vp_{m,h}^k\right\|_{6} \left(h^{1-d/6}\left\|\ub_{m,h}^{k+1}\right\|_{H^1} + \left\|\ub_{m,h}^{k+1}\right\|_{H^1}\right) \left\|\sigma_m^{\ub,k+1}\right\| \non\\
&\leq& C\left\|\vp_m^{k+1}-\vp_{m,h}^k\right\|_{H^1} \left\|\sigma_m^{\ub,k+1}\right\| \non\\
&\leq& C\left\|\vp^{k+1}-\vp_h^k\right\|_{H^1} \left\|\sigma_m^{\ub,k+1}\right\| \non\\
&\leq& \frac{C}{\theta_4}\left( R^{k+1} + \left\|\nabla e^{\vp,k}\right\|^2 \right) + \theta_4\left\|\sigma_m^{\ub,k+1}\right\|^2.
\end{eqnarray}
By choosing $\theta_2 = \frac{\nu_0}{2}, \theta_3 = \alpha_{BJSJ}\frac{\nu_0}{2\sqrt{d\lambda}}, \theta_4 = \frac{\nu_0}{4\lambda}$, we complete the proof of the lemma. \hfill 
\end{proof}

The next lemma contains the estimates of $I_j, j=5, 6, 7, 8$.
\begin{lemma}[Estimates of $I_5, I_6, I_7, I_8$]\label{lemma_RHS5-8}
The assumptions are the same as in Lemma \ref{lemma_RHS1}.  One has the following estimates on  the terms $I_5, I_6, I_7, I_8$ of RHS in \eqref{pd_err}: 
\begin{eqnarray}
\left|\frac{\rho_0}{\chi}\left( R_m^{\ub,k+1},\sigma_m^{\ub,k+1}\right)_m\right| &\leq& C\left\|R_m^{\ub,k+1}\right\|^2 + \frac{\nu_0}{4\lambda}\Big\|\sigma_m^{\ub,k+1}\Big\|^2, \label{RHS5}\\
\left|\rho_0\left( R_c^{\ub,k+1},\sigma_c^{\ub,k+1}\right)_c\right| &\leq& C\left\| R_c^{\ub,k+1}\right\|^2 + \frac{\nu_0}{2}\Big\|\mathbb{D}(\sigma_c^{\ub,k+1})\Big\|^2, \label{RHS6}\\
\left|\left( R^{\varphi,k+1}, \sigma^{\mu,k+1}\right)\right| &\leq& C\left\| R^{\varphi,k+1}\right\|^2 + \frac{M_0}{12} \left\|\nabla\sigma^{\mu,k+1}\right\|^2, \label{RHS7}\\
\left|\left(\rho^{\mu,k+1},\delta_t\sigma^{\varphi,k+1}\right)\right| &\leq& \frac{C}{\theta_8}\Big\|\nabla\rho^{\mu,k+1}\Big\|^2 + \theta_8\Big\|\delta_t\sigma^{\varphi,k+1}\Big\|_{-1,h}^2. \label{RHS8}
\end{eqnarray}
\end{lemma}
\begin{proof}
In fact, \eqref{RHS5} is a direct result of the Cauchy-Schwarz inequality. Thanks to the Poincar\'{e} inequality and Korn's inequality \cite{BrSc2008}, for any $\theta_6>0$, we have
\begin{eqnarray}
\left|\rho_0\left( R_c^{\ub,k+1},\sigma_c^{\ub,k+1}\right)_c\right|
&\leq& \Big\|\rho_0 R_c^{\ub,k+1}\Big\| \Big\|\sigma_c^{\ub,k+1}\Big\| \leq C\Big\|\rho_0 R_c^{\ub,k+1}\Big\| \Big\|\mathbb{D}(\sigma_c^{\ub,k+1})\Big\| \nonumber\\
&\leq& \frac{C}{\theta_6}\Big\| R_c^{\ub,k+1}\Big\|^2 + \theta_6\Big\|\mathbb{D}(\sigma_c^{\ub,k+1})\Big\|^2 ,  \quad \forall \theta_6 > 0 . 
\end{eqnarray}
We notice that $\left( R^{\varphi,k+1}, 1\right) = 0$ holds for all $0\leq k\leq K-1$ by choosing the test function $v=1$ in \eqref{weak_CH1} and using the mass conservation of Ritz projection. Let $\overline{\sigma^{\mu,k+1}}$ be the mean value of $\sigma^{\mu,k+1}$ on $\Omega$, it follows that 
\begin{eqnarray}
\left|\left( R^{\varphi,k+1}, \sigma^{\mu,k+1}\right)\right| &=& \left|\left( R^{\varphi,k+1}, \sigma^{\mu,k+1}-\overline{\sigma^{\mu,k+1}}\right)\right| \leq \Big\| R^{\varphi,k+1}\Big\| \Big\|\sigma^{\mu,k+1}-\overline{\sigma^{\mu,k+1}}\Big\| \nonumber\\
&\leq& C \Big\| R^{\varphi,k+1}\Big\| \Big\|\nabla\sigma^{\mu,k+1}\Big\| \leq \frac{C}{\theta_7}\Big\| R^{\varphi,k+1}\Big\|^2 + \theta_7 \Big\|\nabla\sigma^{\mu,k+1}\Big\|^2 ,  \quad \forall \theta_7 > 0 . 
\end{eqnarray}
For the eighth term of the RHS of \eqref{pd_err}, we  apply Lemma \ref{nabla_-1,h} and recall $\delta_t\sigma^{\varphi,k+1} \in \mathring{Y}_h$ for all $0 \leq k \leq K-1$. Thus  for any $\theta_8>0$, one gets  
\begin{eqnarray}
\left|\left(\rho^{\mu,k+1},\delta_t\sigma^{\varphi,k+1}\right)\right| \leq C\Big\|\nabla\rho^{\mu,k+1}\Big\| \Big\|\delta_t\sigma^{\varphi,k+1}\Big\|_{-1,h} \leq \frac{C}{\theta_8}\Big\|\nabla\rho^{\mu,k+1}\Big\|^2 + \theta_8\Big\|\delta_t\sigma^{\varphi,k+1}\Big\|_{-1,h}^2.
\end{eqnarray}
The proof is complete upon setting $\theta_6 = \frac{\nu_0}{2}, \theta_7 = \frac{M_0}{12}$. \hfill 
\end{proof}

The following lemma gives an estimate of the ninth term $I_9$ on the RHS of \eqref{pd_err}.
\begin{lemma}[Estimate of $I_9$]\label{lemma_RHS9}
The assumptions are the same as in Lemma \ref{lemma_RHS1}. Then for any $0\leq k\leq K-1$, the following inequality holds for a constant $C$ that is independent of $\tau$ and $h$: 
\begin{eqnarray}\label{RHS9}
&&\left(\ub^{k+1}\varphi^{k+1}-\overline{\ub}_h^{k+1}\varphi_h^k, \nabla\sigma^{\mu,k+1}\right) \nonumber\\
&\leq& - \frac{\tau}{\rho_0}\left\|\vp_{c,h}^k\nabla\sigma_c^{\mu,k+1}\right\|^2 - \frac{\tau\chi}{\rho_0}\left\|\vp_{m,h}^k\nabla\sigma_m^{\mu,k+1}\right\|^2 + \frac{M_0}{12} \left\|\nabla\sigma^{\mu,k+1}\right\|^2 + C\Bigg[ \tau^2 \left(1 + \left\|\nabla\rho^{\mu,k+1}\right\|_6^2\right) \non\\
&& + R^{k+1} + \left\|\nabla e^{\varphi,k}\right\|^2 + \left\|\vp_{h}^{k}\right\|_{\infty}^2 \left(R^{k+1} + \left\|e^{\ub,k}\right\|^2\right)\Bigg] . 
\end{eqnarray} 
\end{lemma}
\begin{proof}
We first split the term on $\Omega_m$ as follows:
\begin{eqnarray}\label{RHS8_split}
&& \ub_m^{k+1}\vp_m^{k+1} - \overline{\ub}_{m,h}^{k+1}\vp_{m,h}^k \non\\
&=& \ub_m^{k+1}\vp_m^{k+1} - \left(\ub_{m,h}^k - \frac{\tau\chi}{\rho_0}\vp_{m,h}^k\nabla\mu_{m,h}^{k+1}\right)\vp_{m,h}^k
= \ub_m^{k+1}\vp_m^{k+1} - \ub_{m,h}^k \vp_{m,h}^k + \frac{\tau\chi}{\rho_0}(\vp_{m,h}^k)^2\nabla\mu_{m,h}^{k+1} \non\\
&=& \ub_m^{k+1}\left(\varphi_m^{k+1}-\varphi_m^{k}+e_m^{\varphi,k}\right) + \varphi_{m,h}^{k}\left(\ub_m^{k+1}-\ub_m^{k}+e_m^{\ub,k}\right) - \frac{\tau\chi}{\rho_0}(\vp_{m,h}^k)^2\left(\nabla\rho_m^{\mu,k+1} + \nabla\sigma_m^{\mu,k+1} - \nabla\mu_m^{k+1}\right) \non\\
&=& I_m - \frac{\tau\chi}{\rho_0}\left(\vp_{m,h}^k\right)^2\nabla\sigma_m^{\mu,k+1},
\end{eqnarray}
where
\begin{eqnarray}
I_m \triangleq \ub_m^{k+1}\left(\varphi_m^{k+1}-\varphi_m^{k}+e_m^{\varphi,k}\right) + \varphi_{m,h}^{k}\left(\ub_m^{k+1}-\ub_m^{k}+e_m^{\ub,k}\right) - \frac{\tau\chi}{\rho_0}(\vp_{m,h}^k)^2\left(\nabla\rho_m^{\mu,k+1} + \nabla\mu_m^{k+1}\right).
\end{eqnarray}
In light of  $\|\vp_h^k\|_{H^1}\leq C$ and $\|e^{\vp,k}\|_{H^1}\leq C\|\nabla e^{\vp,k}\|$, one has 
\begin{eqnarray}
\left\|I_m\right\| &\leq& \left\|\ub_m^{k+1}\left(\varphi_m^{k+1}-\varphi_m^{k}+e_m^{\varphi,k}\right)\right\| + \left\|\vp_{m,h}^{k}\left(\ub_m^{k+1}-\ub_m^{k}+e_m^{\ub,k}\right)\right\| + \left\|\frac{\tau\chi}{\rho_0}(\vp_{m,h}^k)^2\left(\nabla\rho_m^{\mu,k+1} + \nabla\mu_m^{k+1}\right)\right\| \non\\
&\leq& \left\|\ub_m^{k+1}\right\|_4 \left\|\varphi_m^{k+1}-\varphi_m^{k}+e_m^{\varphi,k}\right\|_4 + \left\|\vp_h^{k}\right\|_{\infty} \left\|\ub_m^{k+1}-\ub_m^{k}+e_m^{\ub,k}\right\| + C\tau\left\|\vp_{m,h}^k\right\|_{6}^2 \left\|\nabla\rho_m^{\mu,k+1} + \nabla\mu_m^{k+1}\right\|_6 \non\\
&\leq& C \left\|\varphi_m^{k+1}-\varphi_m^{k}+e_m^{\varphi,k}\right\|_{H^1} + \left\|\vp_{h}^{k}\right\|_{\infty} \left(\left\|\ub_m^{k+1}-\ub_m^{k}\right\| + \left\|e_m^{\ub,k}\right\|\right) + C\tau\left\|\vp_{m,h}^k\right\|_{H^1}^2 \left(\left\|\nabla\rho_m^{\mu,k+1}\right\|_6 + C\right) \non\\
&\leq& C \left\|\varphi_m^{k+1}-\varphi_m^{k}\right\|_{H^1} + C\left\|\nabla e_m^{\varphi,k}\right\| + \left\|\vp_{h}^{k}\right\|_{\infty} \left(\left\|\ub_m^{k+1}-\ub_m^{k}\right\| + \left\|e_m^{\ub,k}\right\|\right) + C\tau\left(1+\left\|\nabla\rho_m^{\mu,k+1}\right\|_6\right).
\end{eqnarray}
By Young's inequality, we obtain
\begin{eqnarray}\label{RHS8_Im}
\left\|I_m\right\|^2 &\leq& C\tau^2 \big(1 + \left\|\nabla\rho_m^{\mu,k+1}\right\|_6^2\big) + C \left\|\varphi_m^{k+1}-\varphi_m^{k}\right\|_{H^1}^2 + C\left\|\nabla e_m^{\varphi,k}\right\|^2 + C\left\|\vp_{h}^{k}\right\|_{\infty}^2 \left(\left\|\ub_m^{k+1}-\ub_m^{k}\right\|^2 + \left\|e_m^{\ub,k}\right\|^2\right) \non\\
&\leq& C\tau^2 \big(1 + \left\|\nabla\rho_m^{\mu,k+1}\right\|_6^2\big) + C R^{k+1} + C\left\|\nabla e_m^{\varphi,k}\right\|^2 + C\left\|\vp_{h}^{k}\right\|_{\infty}^2 \left(R^{k+1} + \left\|e_m^{\ub,k}\right\|^2\right).
\end{eqnarray}
Similarly, with the following definition 
\begin{eqnarray}
I_c \triangleq \ub_c^{k+1}\left(\varphi_c^{k+1}-\varphi_c^{k}+e_c^{\varphi,k}\right) + \varphi_{c,h}^{k}\left(\ub_c^{k+1}-\ub_c^{k}+e_c^{\ub,k}\right) - \frac{\tau}{\rho_0}(\vp_{c,h}^k)^2\left(\nabla\rho_c^{\mu,k+1} + \nabla\mu_c^{k+1}\right),
\end{eqnarray}
one gets 
\begin{eqnarray}\label{RHS8_Ic}
\left\|I_c\right\|^2 \leq C\tau^2 \big(1 + \left\|\nabla\rho_c^{\mu,k+1}\right\|_6^2\big) + C R^{k+1} + C\left\|\nabla e_c^{\varphi,k}\right\|^2 + C\left\|\vp_{h}^{k}\right\|_{\infty}^2 \left(R^{k+1} + \left\|e_c^{\ub,k}\right\|^2\right).
\end{eqnarray}
Consequently, the following inequality is valid: 
\begin{eqnarray}\label{RHS8_IcIm}
\left\|I_c\right\|^2 + \left\|I_m\right\|^2 \leq C\tau^2 \big(1 + \left\|\nabla\rho^{\mu,k+1}\right\|_6^2\big) + C R^{k+1} + C\left\|\nabla e^{\varphi,k}\right\|^2 + C\left\|\vp_{h}^{k}\right\|_{\infty}^2 \left(R^{k+1} + \left\|e^{\ub,k}\right\|^2\right).
\end{eqnarray}
Thus for constant $\theta_9>0$, there holds
\begin{eqnarray}
&& \left(\ub^{k+1}\varphi^{k+1}-\overline{\ub}_h^{k+1}\varphi_h^k, \nabla\sigma^{\mu,k+1}\right) \non\\
&=& \left(\ub_c^{k+1}\varphi_c^{k+1}-\overline{\ub}_{c,h}^{k+1}\varphi_{c,h}^k, \nabla\sigma_c^{\mu,k+1}\right)_c + \left(\ub_m^{k+1}\varphi_m^{k+1}-\overline{\ub}_{m,h}^{k+1}\varphi_{m,h}^k, \nabla\sigma_m^{\mu,k+1}\right)_m \non\\
&=& \left(I_c - \frac{\tau}{\rho_0}\left(\vp_{c,h}^k\right)^2\nabla\sigma_c^{\mu,k+1}, \nabla\sigma_c^{\mu,k+1}\right)_c + \left(I_m - \frac{\tau\chi}{\rho_0}\left(\vp_{m,h}^k\right)^2\nabla\sigma_m^{\mu,k+1}, \nabla\sigma_m^{\mu,k+1}\right)_m \non\\
&=& \left(I_c, \nabla\sigma_c^{\mu,k+1}\right)_c + \left(I_m, \nabla\sigma_m^{\mu,k+1}\right)_m - \frac{\tau}{\rho_0}\left\|\vp_{c,h}^k\nabla\sigma_c^{\mu,k+1}\right\|^2 - \frac{\tau\chi}{\rho_0}\left\|\vp_{m,h}^k\nabla\sigma_m^{\mu,k+1}\right\|^2 \non\\
&\leq& - \frac{\tau}{\rho_0}\left\|\vp_{c,h}^k\nabla\sigma_c^{\mu,k+1}\right\|^2 - \frac{\tau\chi}{\rho_0}\left\|\vp_{m,h}^k\nabla\sigma_m^{\mu,k+1}\right\|^2 + \frac{1}{4\theta_9}\left(\left\|I_c\right\|^2 + \left\|I_m\right\|^2\right) + \theta_9 \left\|\nabla\sigma^{\mu,k+1}\right\|^2 \non\\
&\leq& - \frac{\tau}{\rho_0}\left\|\vp_{c,h}^k\nabla\sigma_c^{\mu,k+1}\right\|^2 - \frac{\tau\chi}{\rho_0}\left\|\vp_{m,h}^k\nabla\sigma_m^{\mu,k+1}\right\|^2 + \theta_9 \left\|\nabla\sigma^{\mu,k+1}\right\|^2 + \frac{C}{\theta_9}\Bigg[ \tau^2 \big(1 + \left\|\nabla\rho^{\mu,k+1}\right\|_6^2\big) \non\\
&& + R^{k+1} + \left\|\nabla e^{\varphi,k}\right\|^2 + \left\|\vp_{h}^{k}\right\|_{\infty}^2 \left(R^{k+1} + \left\|e^{\ub,k}\right\|^2\right)\Bigg].
\end{eqnarray}
This proves the lemma by choosing $\theta_9 = \frac{M_0}{12}$. \hfill 
\end{proof}

The term $I_{10}$ is estimated in the following lemma.
\begin{lemma}[Estimate of the  term $I_{10}$]\label{lemma_RHS10}
The assumptions are the same as in Lemma \ref{lemma_RHS1}.. Then the tenth term of RHS in \eqref{pd_err} satisfies
\begin{eqnarray}\label{RHS10}
&& \left|\frac{\gamma}{\epsilon}\left(f(\varphi^{k+1},\varphi^{k+1})-f(\varphi_h^{k+1},\varphi_h^k),
\delta_t\sigma^{\varphi,k+1}\right)\right| \non\\
&\leq& \theta_{10} \left\|\delta_t\sigma^{\varphi,k+1}\right\|_{-1,h}^2 + \frac{C}{\theta_{10}}\left(R^{k+1} + \left(1+\left\|\varphi_h^{k+1}\right\|_{\infty}^{4}\right)\left\|\nabla e^{\varphi,k+1}\right\|^2 + \left\|\nabla e^{\varphi,k}\right\|^2\right) , 
\end{eqnarray}
for a constant $C$  independent of $\tau$ and $h$.
\end{lemma}
\begin{proof}
First we need to estimate $\left\|\nabla \left(f(\varphi^{k+1},\varphi^{k+1})-f(\varphi_h^{k+1},\varphi_h^k)\right)\right\|$. Recall that $f(a,b)=a^3-b$. Hence 
\begin{eqnarray}
&& \Big\|\nabla \left(f(\varphi^{k+1},\varphi^{k+1})-f(\varphi_h^{k+1},\varphi_h^k)\right)\Big\| \nonumber\\
&\leq& \Big\|\nabla(\varphi^{k+1})^{3} - \nabla(\varphi_h^{k+1})^{3}\Big\| + \Big\|\nabla(\varphi^{k+1}-\varphi_h^{k})\Big\| \nonumber\\
&=& \Big\|3(\varphi^{k+1})^{2}\nabla\varphi^{k+1} - 3(\varphi_h^{k+1})^{2}\nabla\varphi_h^{k+1}\Big\| + \Big\|\nabla\left(\varphi^{k+1}-\varphi^{k}+\varphi^{k}-\varphi_h^{k}\right)\Big\| \nonumber\\
&\leq& 3\Big\|\left((\varphi^{k+1})^{2}-(\varphi_h^{k+1})^{2}\right)\nabla\varphi^{k+1} + (\varphi_h^{k+1})^{2}\nabla\left(\varphi^{k+1}-\varphi_h^{k+1}\right)\Big\| + \Big\|\nabla\left(\vp^{k+1}-\vp^k\right)\Big\| + \Big\|\nabla e^{\varphi,k}\Big\| \nonumber\\
&\leq& 3\Big\|\varphi^{k+1}+\varphi_h^{k+1}\Big\|_{6}\Big\| e^{\varphi,k+1}\Big\|_{6}\Big\|\nabla\varphi^{k+1}\Big\|_{6} + 3\Big\|\varphi_h^{k+1}\Big\|_{\infty}^{2}\Big\|\nabla e^{\varphi,k+1}\Big\| + \Big\|\nabla\left(\vp^{k+1}-\vp^k\right)\Big\| + \Big\|\nabla e^{\varphi,k}\Big\| \nonumber\\
&\leq& C\left(\left\|\varphi^{k+1}\right\|_6+C\left\|\varphi_h^{k+1}\right\|_{H_1}\right)\Big\|\nabla e^{\varphi,k+1}\Big\| + 3\Big\|\varphi_h^{k+1}\Big\|_{\infty}^{2}\Big\|\nabla e^{\varphi,k+1}\Big\| + \Big\|\nabla\left(\vp^{k+1}-\vp^k\right)\Big\| + \Big\|\nabla e^{\varphi,k}\Big\| \nonumber\\
&\leq& \Big\|\nabla\left(\vp^{k+1}-\vp^k\right)\Big\| + C\left(1+\Big\|\varphi_h^{k+1}\Big\|_{\infty}^{2}\right)\Big\|\nabla e^{\varphi,k+1}\Big\| + \Big\|\nabla e^{\varphi,k}\Big\|,
\end{eqnarray}
which in turn yields
\begin{eqnarray}
&& \Big\|\nabla \left(f(\varphi^{k+1},\varphi^{k+1})-f(\varphi_h^{k+1},\varphi_h^k)\right)\Big\|^2 \nonumber\\
&\leq& C\Big\|\nabla\left(\vp^{k+1}-\vp^k\right)\Big\|^2 + C\left(1+\Big\|\varphi_h^{k+1}\Big\|_{\infty}^{2}\right)^2\Big\|\nabla e^{\varphi,k+1}\Big\|^2 + C\Big\|\nabla e^{\varphi,k}\Big\|^2 \nonumber\\
&\leq& CR^{k+1} + C\left(1+\left\|\varphi_h^{k+1}\right\|_{\infty}^{4}\right)\Big\|\nabla e^{\varphi,k+1}\Big\|^2 + C\Big\|\nabla e^{\varphi,k}\Big\|^2.
\end{eqnarray}
Thus by Lemma \ref{nabla_-1,h}, we derive the following estimate for any $\theta_{10} >0$: 
\begin{eqnarray}
&& \left|\frac{\gamma}{\epsilon}\left(f(\varphi^{k+1},\varphi^{k+1})-f(\varphi_h^{k+1},\varphi_h^k),
\delta_t\sigma^{\varphi,k+1}\right)\right| \nonumber\\
&\leq& C\left\|\nabla \left(f(\varphi^{k+1},\varphi^{k+1})-f(\varphi_h^{k+1},\varphi_h^k)\right)\right\| \left\|\delta_t\sigma^{\varphi,k+1}\right\|_{-1,h} \nonumber\\
&\leq& \theta_{10}\left\|\delta_t\sigma^{\varphi,k+1}\right\|_{-1,h}^2 + \frac{C}{\theta_{10}} \left\|\nabla\left(f(\varphi^{k+1},\varphi^{k+1})-f(\varphi_h^{k+1},\varphi_h^k)\right)\right\|^2 \nonumber\\
&\leq& \theta_{10} \left\|\delta_t\sigma^{\varphi,k+1}\right\|_{-1,h}^2 + \frac{C}{\theta_{10}}\left(R^{k+1} + \left(1+\left\|\varphi_h^{k+1}\right\|_{\infty}^{4}\right)\left\|\nabla e^{\varphi,k+1}\right\|^2 + \left\|\nabla e^{\varphi,k}\right\|^2\right).
\end{eqnarray}
This completes the proof.
\hfill 
\end{proof}

Finally we estimate the last term $I_{11}$ in the following lemma.
\begin{lemma}[Estimate of the $I_{11}$]\label{lemma_RHS11}
The assumptions are the same as in Lemma \ref{lemma_RHS1}. Then for the last term $I_{11}$ of RHS in \eqref{pd_err}, 
the following inequality holds for a constant $C$ independent of $\tau$ and $h$: 
\begin{eqnarray}\label{RHS11}
&& \left| - \left(\varphi^{k+1}\nabla\mu^{k+1}-\varphi_{h}^k\nabla\mu_{h}^{k+1},\sigma^{\ub,k+1}\right) \right| \non\\
&\leq& C\left(R^{k+1} + \left\|\nabla e^{\vp,k}\right\|^2\right) + \frac{M_0}{12}\left\|\nabla e^{\mu,k+1}\right\|^2 + \left(1 + C\left\|\varphi_h^{k}\right\|_{\infty}^2\right) \left\|\sigma^{\ub,k+1}\right\|^2.
\end{eqnarray}
\end{lemma}
\begin{proof}
We make use of the following decomposition
\begin{eqnarray}
\left\|\varphi^{k+1}\nabla\mu^{k+1} - \varphi_h^{k}\nabla\mu_h^{k+1}\right\| &=& \left\|(\varphi^{k+1}-\varphi_h^{k})\nabla\mu^{k+1} + \varphi_h^{k}\nabla(\mu^{k+1}-\mu_h^{k+1})\right\| \nonumber\\
&=& \left\|\left(\varphi^{k+1}-\varphi^{k}+e^{\varphi,k}\right)\nabla\mu^{k+1} + \varphi_h^{k}\nabla e^{\mu,k+1}\right\| \nonumber\\
&\leq& \left\|\varphi^{k+1}-\varphi^{k}+e^{\varphi,k}\right\|_4 \left\|\nabla\mu^{k+1}\right\|_{4} + \left\|\varphi_h^{k}\right\|_{\infty}\left\|\nabla e^{\mu,k+1}\right\| \nonumber\\
&\leq& C\left(\left\|\vp^{k+1}-\vp^k\right\|_4 + \left\|e^{\varphi,k}\right\|_4\right) + \Big\|\varphi_h^{k}\Big\|_{\infty} \left\|\nabla e^{\mu,k+1}\right\| \non\\
&\leq& C\left(\left\|\vp^{k+1}-\vp^k\right\|_{H^1} + \left\|e^{\vp,k}\right\|_{H^1}\right) + \left\|\vp_h^k\right\|_{\infty} \left\|\nabla e^{\mu,k+1}\right\|.
\end{eqnarray}

Then for any $\theta_{11}>0$ there holds
\begin{eqnarray}
&& \left|\left(\varphi^{k+1}\nabla\mu^{k+1}-\varphi_h^k\nabla\mu_h^{k+1},\sigma^{\ub,k+1}\right)\right| \leq \Big\|\varphi^{k+1}\nabla\mu^{k+1}-\varphi_h^k\nabla\mu_h^{k+1}\Big\| \Big\|\sigma^{\ub,k+1}\Big\| \non\\
&\leq& \Bigg[C\left(\left\|\vp^{k+1}-\vp^k\right\|_{H^1} + \left\|e^{\vp,k}\right\|_{H^1}\right) + \left\|\vp_h^k\right\|_{\infty} \left\|\nabla e^{\mu,k+1}\right\|\Bigg] \left\|\sigma^{\ub,k+1}\right\| \non\\
&\leq& C\left(\left\|\vp^{k+1}-\vp^k\right\|_{H^1} + \left\|e^{\vp,k}\right\|_{H^1}\right)\left\|\sigma^{\ub,k+1}\right\| + \left\|\nabla e^{\mu,k+1}\right\| \left\|\varphi_h^{k}\right\|_{\infty}\left\|\sigma^{\ub,k+1}\right\| \non\\
&\leq& C\left(\left\|\vp^{k+1}-\vp^k\right\|_{H^1}^2 + \left\|e^{\vp,k}\right\|_{H^1}^2\right) + \left\|\sigma^{\ub,k+1}\right\|^2 +\theta_{11}\left\|\nabla e^{\mu,k+1}\right\|^2 + \frac{C}{\theta_{11}}\Big\|\varphi_h^{k}\Big\|_{\infty}^2 \Big\|\sigma^{\ub,k+1}\Big\|^2 \non\\
&\leq& C\left(R^{k+1} + \left\|\nabla e^{\vp,k}\right\|^2\right) + \theta_{11}\left\|\nabla e^{\mu,k+1}\right\|^2 + \left(1 + \frac{C}{\theta_{11}}\left\|\varphi_h^{k}\right\|_{\infty}^2\right) \left\|\sigma^{\ub,k+1}\right\|^2.
\end{eqnarray}
The proof is complete by choosing $\theta_{11} = \frac{M_0}{12}$. \hfill 
\end{proof}

The next lemma gives an estimate of $\left\|\delta_t\sigma^{\varphi,k+1}\right\|_{-1,h}$.
\begin{lemma}\label{-1,h_norm_estimate}
The assumptions are the same as in Lemma \ref{lemma_RHS1}. There exists a constant $C>0$  independent of $\tau$ and $h$ such that
\begin{eqnarray}\label{negative-es}
\left\|\delta_t\sigma^{\varphi,k+1}\right\|_{-1,h}^2 &\leq& C\tau^2 + C\tau^2\left\|\nabla\rho^{\mu,k+1}\right\|_6^2 + \left(\frac{25M_1^2}{4} + C_1\tau(T+1)\right) \Big\|\nabla\sigma^{\mu,k+1}\Big\|^2  + C\Big\| R^{\varphi,k+1}\Big\|^2 \non\\
&& + C\left(1+\left\|\vp_h^k\right\|_{\infty}^2\right)R^{k+1} + C\left\|\nabla e^{\vp,k}\right\|^2 + C\left\|\vp_h^k\right\|_{\infty}^2\left\|e^{\ub,k}\right\|^2.
\end{eqnarray}
\end{lemma}

\begin{proof}
Recall that $\|\zeta\|_{-1,h}^2 = \|\nabla \textsf{T}_h(\zeta)\|^2 = \big(\nabla \textsf{T}_h(\zeta), \nabla \textsf{T}_h(\zeta)\big) = \big(\zeta, \textsf{T}_h(\zeta)\big)$ for all $\zeta\in\mathring{Y}_h$. Noticing that $\delta_t\sigma^{\varphi,k+1} \in \mathring{Y}_h$, setting $v=\textsf{T}_h\big(\delta_t\sigma^{\varphi,k+1}\big)$ in \eqref{pd_err_CH1}, using \eqref{RHS1}, \eqref{RHS8_split} and \eqref{RHS8_IcIm}, we derive
\begin{eqnarray}
&& \left\|\delta_t\sigma^{\varphi,k+1}\right\|_{-1,h}^2 = \left(\delta_t\sigma^{\varphi,k+1}, \textsf{T}_h\big(\delta_t\sigma^{\varphi,k+1}\big)\right) \nonumber\\
&=& - \left(\left({\rm M}(\vp^{k+1})-{\rm M}(\vp_h^k)\right)\nabla\widetilde{\mathcal{P}}^{k+1}\mu^{k+1}, \nabla\textsf{T}_h\big(\delta_t\sigma^{\varphi,k+1}\big)\right) + \left( R^{\varphi,k+1}, \textsf{T}_h\big(\delta_t\sigma^{\varphi,k+1}\big)\right) \non\\
&& -\left({\rm M}(\vp_h^k)\nabla\sigma^{\mu,k+1},\nabla \textsf{T}_h\big(\delta_t\sigma^{\varphi,k+1}\big)\right) + \left(\ub^{k+1}\varphi^{k+1}-\overline{\ub}_h^{k+1}\varphi_h^k,\nabla \textsf{T}_h\big(\delta_t\sigma^{\varphi,k+1}\big)\right) \nonumber\\
&\leq& C\left(R^{k+1} + \left\|\nabla e^{\vp,k}\right\|^2\right) + \frac{1}{5}\left\|\nabla\textsf{T}_h\big(\delta_t\sigma^{\varphi,k+1}\big)\right\|^2 + \left\| R^{\varphi,k+1}\right\| \left\|\textsf{T}_h\big(\delta_t\sigma^{\varphi,k+1}\big)\right\| \non\\
&& + \left\|{\rm M}(\vp_h^k)\nabla\sigma^{\mu,k+1}\right\| \left\|\nabla\textsf{T}_h\big(\delta_t\sigma^{\varphi,k+1}\big)\right\| + \left\|\ub^{k+1}\varphi^{k+1}-\overline{\ub}_h^{k+1}\varphi_h^k\right\| \left\|\nabla\textsf{T}_h\big(\delta_t\sigma^{\varphi,k+1}\big)\right\| \nonumber\\
&\leq& CR^{k+1} + C\left\|\nabla e^{\vp,k}\right\|^2 + \frac{1}{5}\left\|\nabla\textsf{T}_h\big(\delta_t\sigma^{\varphi,k+1}\big)\right\|^2 + C\left\| R^{\varphi,k+1}\right\| \left\|\nabla\textsf{T}_h\big(\delta_t\sigma^{\varphi,k+1}\big)\right\| \non\\
&& + \frac{5}{4}\left\|{\rm M}(\vp_h^k)\nabla\sigma^{\mu,k+1}\right\|^2 + \frac{1}{5}\left\|\nabla\textsf{T}_h\big(\delta_t\sigma^{\varphi,k+1}\big)\right\|^2 + \frac{5}{4}\left\|\ub^{k+1}\varphi^{k+1}-\overline{\ub}_h^{k+1}\varphi_h^k\right\|^2 + \frac{1}{5}\left\|\nabla\textsf{T}_h\big(\delta_t\sigma^{\varphi,k+1}\big)\right\|^2 \nonumber\\
&\leq& CR^{k+1} + C\left\|\nabla e^{\vp,k}\right\|^2 + \frac{5M_1^2}{4}\left\|\nabla\sigma^{\mu,k+1}\right\|^2 + C\left\| R^{\varphi,k+1}\right\|^2 + \frac{4}{5}\left\|\nabla\textsf{T}_h\big(\delta_t\sigma^{\varphi,k+1}\big)\right\|^2  \non\\
&& + \frac{5}{4}\left\|I_c - \frac{\tau}{\rho_0}(\varphi_{c,h}^{k})^2\nabla\sigma_c^{\mu,k+1}\right\|^2 + \frac{5}{4}\left\|I_m - \frac{\tau\chi}{\rho_0}(\varphi_{m,h}^{k})^2\nabla\sigma_m^{\mu,k+1}\right\|^2 \nonumber\\
&\leq& CR^{k+1} + C\left\|\nabla e^{\vp,k}\right\|^2 + \frac{5M_1^2}{4}\left\|\nabla\sigma^{\mu,k+1}\right\|^2 + C\left\|R^{\varphi,k+1}\right\|^2 + \frac{4}{5}\Big\|\delta_t\sigma^{\varphi,k+1}\Big\|_{-1,h}^2 \non\\
&& + \frac{5}{2}\left\|I_c\right\|^2 + \frac{5}{2}\left\| \frac{\tau}{\rho_0}(\varphi_{c,h}^{k})^2\nabla\sigma_c^{\mu,k+1}\right\|^2 + \frac{5}{2}\left\|I_m\right\|^2 + \frac{5}{2}\left\|\frac{\tau\chi}{\rho_0}(\varphi_{m,h}^{k})^2\nabla\sigma_m^{\mu,k+1}\right\|^2 \nonumber\\
&\leq& CR^{k+1} + C\left\|\nabla e^{\vp,k}\right\|^2 + \frac{5M_1^2}{4}\left\|\nabla\sigma^{\mu,k+1}\right\|^2 + C\left\| R^{\varphi,k+1}\right\|^2 + \frac{4}{5}\left\|\delta_t\sigma^{\varphi,k+1}\right\|_{-1,h}^2 \non\\
&& + \frac{5}{2}\left(\left\|I_c\right\|^2 + \left\|I_m\right\|^2\right) + \frac{5\tau^2}{2\rho_0^2}\left\|\varphi_{h}^{k}\right\|_{\infty}^4 \left\|\nabla\sigma_c^{\mu,k+1}\right\|^2 + \frac{5\tau^2\chi^2}{2\rho_0^2}\left\|\varphi_{h}^{k}\right\|_{\infty}^4 \left\|\nabla\sigma_m^{\mu,k+1}\right\|^2 \nonumber\\
&\leq& CR^{k+1} + C\left\|\nabla e^{\vp,k}\right\|^2 + \left(\frac{5M_1^2}{4} + C\tau^2\left\|\varphi_{h}^{k}\right\|_{\infty}^4\right)\left\|\nabla\sigma^{\mu,k+1}\right\|^2 + C\Big\| R^{\varphi,k+1}\Big\|^2 + \frac{4}{5}\Big\|\delta_t\sigma^{\varphi,k+1}\Big\|_{-1,h}^2 \non\\
&& + C\tau^2 + C\tau^2\left\|\nabla\rho^{\mu,k+1}\right\|_6^2 + C\left(1+\left\|\vp_h^k\right\|_{\infty}^2\right)R^{k+1} + C\left\|\nabla e^{\vp,k}\right\|^2 + C\left\|\vp_h^k\right\|_{\infty}^2 \left\|e^{\ub,k}\right\|^2 \non\\
&\leq& C\tau^2 + \left(\frac{5M_1^2}{4} + C\tau^2\left\|\varphi_{h}^{k}\right\|_{\infty}^4\right) \Big\|\nabla\sigma^{\mu,k+1}\Big\|^2  + C\Big\| R^{\varphi,k+1}\Big\|^2 + \frac{4}{5}\left\|\delta_t\sigma^{\varphi,k+1}\right\|_{-1,h}^2 \non\\
&& + C\left(1+\left\|\vp_h^k\right\|_{\infty}^2\right)R^{k+1} + C\left\|\nabla e^{\vp,k}\right\|^2 + C\tau^2\left\|\nabla\rho^{\mu,k+1}\right\|_6^2 + C\left\|\vp_h^k\right\|_{\infty}^2\left\|e^{\ub,k}\right\|^2.
\end{eqnarray}
Since $\tau\left\|\vp_h^k\right\|_{\infty}^4\leq \tau + \tau\left\|\vp_h^k\right\|_{\infty}^\frac{4(6-d)}{d}\leq C(T+1)$  from Lemma \ref{numerical_estimate},
the proof  is complete once we  move $\frac{4}{5}\left\|\delta_t\sigma^{\varphi,k+1}\right\|_{-1,h}^2$ to the left-hand side of the inequality. 
\end{proof}

With all these estimates of the RHS terms in place, the error equation \eqref{pd_err} leads to the following result. 
\begin{lemma}\label{first_combine_lamma}
Suppose $(\varphi, \mu, \ub_c, \ub_m, P_c, P_m)$ is a weak solution to \eqref{ww_CH1}--\eqref{ww_D} satisfying additional regularities prescribed in Assumption \ref{higher_regularities}. Then, for any $\tau, h>0$,  there exists a constant $C>0$, independent of $h$ and $\tau$, such that for any $0\leq k \leq K-1$,
\begin{eqnarray}\label{combine}
&& \frac{M_0}{3}\left\|\nabla\sigma^{\mu,k+1}\right\|^2 + \frac{\gamma\epsilon}{2\tau}\left(\left\|\nabla\sigma^{\varphi,k+1}\right\|^2 - \left\|\nabla\sigma^{\varphi,k}\right\|^2 + \left\|\nabla(\sigma^{\varphi,k+1}-\sigma^{\varphi,k})\right\|^2\right)\nonumber\\
&&\ \ +\ \frac{\rho_0}{2\tau}\left(\left\|\sigma_c^{\ub,k+1}\right\|^2 - \left\|\sigma_c^{\ub,k}\right\|^2 + \left\|\sigma_c^{\ub,k+1}-\sigma_c^{\ub,k}\right\|^2\right)\nonumber\\
&&\ \ +\ \nu_0\left\|\mathbb{D}(\sigma_c^{\ub,k+1})\right\|^2 + \alpha_{BJSJ}\frac{\nu_0}{2\sqrt{d\lambda}}\sum_{i=1}^{d-1}\left\|\sigma_c^{\ub,k+1}\cdot\btau_i\right\|_{cm}^2\nonumber\\
&&\ \ +\ \frac{\rho_0}{2\tau\chi}\left(\left\|\sigma_m^{\ub,k+1}\right\|^2 - \left\|\sigma_m^{\ub,k}\right\|^2 + \left\|\sigma_m^{\ub,k+1}-\sigma_m^{\ub,k}\right\|^2\right) + \frac{\nu_0}{2\lambda}\left\|\sigma_m^{\ub,k+1}\right\|^2 \nonumber\\
&&\ \ +\ \frac{\tau}{\rho_0}\Big\|\varphi_{c,h}^{k}\nabla\sigma_c^{\mu,k+1}\Big\|^2 + \frac{\tau\chi}{\rho_0}\Big\|\varphi_{m,h}^{k}\nabla\sigma_m^{\mu,k+1}\Big\|^2 \non\\
&\leq& \ C\mathcal{R}^{k+1} + C\left(1+\left\|\vp_h^{k+1}\right\|_{\infty}^{4}\right)\left\|\nabla\sigma^{\vp,k+1}\right\|^2 + \left(1+C\left\|\varphi_h^k\right\|_{\infty}^2\right) \left\|\sigma^{\ub,k+1}\right\|^2 \non\\
&& \ + C\left\|\vp_h^k\right\|_{\infty}^2 \left\|\sigma^{\ub,k}\right\|^2 + C\left\|\nabla\sigma^{\vp,k}\right\|^2,
\end{eqnarray}
where
\begin{eqnarray}\label{R_eq}
\mathcal{R}^{k+1} &:=& \tau^2 + \left\|R_m^{\ub,k+1}\right\|^2 + \left\|R_c^{\ub,k+1}\right\|^2 + \left\|R^{\vp,k+1}\right\|^2 + \left(1+\left\|\vp_h^k\right\|_{\infty}^2\right)R^{k+1} \non\\
&& +\ \left\|\vp_h^k\right\|_{\infty}^2 \left\|\rho^{\ub,k}\right\|^2 + \left\|\nabla\rho^{\vp,k}\right\|^2 + \left(1+\left\|\vp_h^{k+1}\right\|_{\infty}^{4}\right)\left\|\nabla\rho^{\vp,k+1}\right\|^2 \non\\
&& +\ (1+\tau^2)\left\|\nabla\rho^{\mu,k+1}\right\|_6^2.
\end{eqnarray}
\end{lemma}
\begin{proof}
Substituting the estimates in Lemmas \ref{lemma_RHS1} - \ref{-1,h_norm_estimate} into the right-hand side of  the error equation \eqref{pd_err}, choosing
\begin{eqnarray}
\theta_8 = \theta_{10} = \frac{M_0}{6\left(\frac{25 M_1^2}{4} + C_1\tau(T+1)\right)},
\end{eqnarray}
with $C_1$ the positive constant defined in inequality \eqref{negative-es}, we get 
\begin{eqnarray}
&& \frac{M_0}{3}\left\|\nabla\sigma^{\mu,k+1}\right\|^2 + \frac{\gamma\epsilon}{2\tau}\left(\|\nabla\sigma^{\varphi,k+1}\|^2 - \|\nabla\sigma^{\varphi,k}\|^2 + \|\nabla(\sigma^{\varphi,k+1}-\sigma^{\varphi,k})\|^2\right)\nonumber\\
&&\ \ +\ \frac{\rho_0}{2\tau}\left(\|\sigma_c^{\ub,k+1}\|^2 - \|\sigma_c^{\ub,k}\|^2 + \|\sigma_c^{\ub,k+1}-\sigma_c^{\ub,k}\|^2\right)\nonumber\\
&&\ \ +\ \nu_0\left\|\mathbb{D}(\sigma_c^{\ub,k+1})\right\|^2 + \alpha_{BJSJ}\frac{\nu_0}{2\sqrt{d\lambda}}\sum_{i=1}^{d-1}\left\|\sigma_c^{\ub,k+1}\cdot\btau_i\right\|_{cm}^2\nonumber\\
&&\ \ +\ \frac{\rho_0}{2\tau\chi}\left(\|\sigma_m^{\ub,k+1}\|^2 - \|\sigma_m^{\ub,k}\|^2 + \|\sigma_m^{\ub,k+1}-\sigma_m^{\ub,k}\|^2\right) + \frac{\nu_0}{2\lambda}\left\|\sigma_m^{\ub,k+1}\right\|^2 \nonumber\\
&&\ \ +\ \frac{\tau}{\rho_0}\Big\|\varphi_{c,h}^{k}\nabla\sigma_c^{\mu,k+1}\Big\|^2 + \frac{\tau\chi}{\rho_0}\Big\|\varphi_{m,h}^{k}\nabla\sigma_m^{\mu,k+1}\Big\|^2 \non\\
&\leq& \ C\tau^2 + C\left\|R_m^{\ub,k+1}\right\|^2 + C\left\|R_c^{\ub,k+1}\right\|^2 + C\left\|R^{\vp,k+1}\right\|^2 + C\left(1+\left\|\vp_h^k\right\|_{\infty}^2\right)R^{k+1} \non\\
&& +\ C\left\|\vp_h^k\right\|_{\infty}^2 \left\|e^{\ub,k}\right\|^2 + C\left\|\nabla e^{\vp,k}\right\|^2 + C\left(1+\left\|\vp_h^{k+1}\right\|_{\infty}^{4}\right)\left\|\nabla e^{\vp,k+1}\right\|^2 \non\\
&& +\ C\left\|\nabla\rho^{\mu,k+1}\right\|^2 + C\tau^2\left\|\nabla\rho^{\mu,k+1}\right\|_6^2 + \left(1+C\left\|\varphi_h^k\right\|_{\infty}^2\right) \left\|\sigma^{\ub,k+1}\right\|^2.
\end{eqnarray}
The proof is complete since  $\|e^{\ub,k}\|^2 = \|\rho^{\ub,k} + \sigma^{\ub,k}\|^2 \leq 2\left(\|\rho^{\ub,k}\|^2 + \|\sigma^{\ub,k}\|^2\right)$,  and $\|\nabla\rho^{\mu,k+1}\| \leq C\|\nabla\rho^{\mu,k+1}\|_6$. \hfill 
\end{proof}

Regarding $\mathcal{R}^{k+1}$ in Eq. \eqref{R_eq}, the following estimate could be derived. 
\begin{lemma}\label{sum-R}
Suppose $(\varphi, \mu, \ub_c, \ub_m, P_c, P_m)$ is a weak solution to \eqref{ww_CH1}--\eqref{ww_D} satisfying additional regularities  in Assumption \ref{higher_regularities}. Then for all $0\leq l\leq K-1$ there holds
\begin{eqnarray}
\sum_{k=0}^{l}\mathcal{R}^{k+1} &\leq& C(T+1)\tau + \frac{2}{\tau}\int_{0}^{t_{l+1}}\left(\left\|\partial_t\rho^{\vp}(\cdot,t)\right\|^2 + \left\|\partial_t\rho^{\ub}(\cdot,t)\right\|^2\right)dt \non\\
&&\ +\ C\tau^{-1/2}\sqrt{T+1}\left(\sum_{k=0}^{l}\left\|\nabla\rho^{\vp,k+1}\right\|^4\right)^{1/2} + \sum_{k=0}^l\left(\left\|\nabla\rho^{\vp,k}\right\|^2 + (1+\tau^2)\left\|\nabla\rho^{\mu,k+1}\right\|_6^2 \right) \non\\
&&\ +\ C\tau^{-1/2}\sqrt{T+1} \left(\sum_{k=0}^{l}\left\|\rho^{\ub,k}\right\|^4\right)^{1/2}.
\end{eqnarray}
\end{lemma}

\begin{proof}
First, by Minkowski's inequality and H\"{o}lder's inequality one obtains
\begin{eqnarray}
&& \left\|R^{\vp,k+1}\right\|^2 = \left\|\delta_t\mathcal{P}\vp^{k+1}-\partial_t\vp^{k+1}\right\|^2 \non\\
&\leq& 2\left\|\delta_t\left(\mathcal{P}\vp^{k+1} - \vp^{k+1}\right)\right\|^2 + 2\left\|\delta_t \vp^{k+1}-\partial_t\vp^{k+1}\right\|^2 \non\\
&=& { \frac{2}{\tau^2}\left\|\int_{t_k}^{t_{k+1}}\partial_t\rho^{\vp}(\cdot,t)dt\right\|^2 } + \frac{2}{\tau^2}\left\|\int_{t_k}^{t_{k+1}}(t-t_k)\partial_{tt}\vp(\cdot,t) dt\right\|^2 \non\\
&\leq& \frac{2}{\tau^2}\left(\int_{t_k}^{t_{k+1}}\left\|\partial_t\rho^{\vp}(\cdot,t)\right\|dt\right)^2 + \frac{2}{\tau^2}\left(\int_{t_k}^{t_{k+1}}(t-t_k)\left\|\partial_{tt}\vp(\cdot,t)\right\| dt\right)^2 \non\\
&\leq& \frac{2}{\tau}\int_{t_k}^{t_{k+1}}\left\|\partial_t\rho^{\vp}(\cdot,t)\right\|^2dt + \frac{2\tau}{3}\int_{t_k}^{t_{k+1}}\left\|\partial_{tt}\vp(\cdot,t)\right\|^2 dt.
\end{eqnarray}
Likewise, for $j\in\{c,m\}$, one has
\begin{eqnarray}
\left\|R_j^{\ub,k+1}\right\|^2 \leq { \frac{2}{\tau}\int_{t_k}^{t_{k+1}}\left\|\partial_t\rho_j^{\ub}(\cdot,t)\right\|^2dt } + \frac{2\tau}{3}\int_{t_k}^{t_{k+1}}\left\|\partial_{tt}\ub(\cdot,t)\right\|^2 dt.
\end{eqnarray}

Applying Minkowski's inequality and H\"{o}lder's inequality again gives, for $j\in\{c,m\}$, 
\begin{eqnarray}
&& \left\|\varphi^{k+1}-\varphi^{k}\right\|^4 = \Big\|\int_{t_k}^{t_{k+1}}\partial_t \varphi(\cdot,t)dt\Big\|^4 \leq \left(\int_{t_k}^{t_{k+1}}\Big\|\partial_t \varphi(\cdot,t)\Big\|dt\right)^4 \non\\
&\leq& \left(\int_{t_k}^{t_{k+1}}\left\|\partial_t\vp(\cdot,t)\right\|^4 dt\right) \left(\int_{t_k}^{t_{k+1}}dt\right)^3 = \tau^3 \int_{t_k}^{t_{k+1}}\left\|\partial_t\vp(\cdot,t)\right\|^4 dt,
\end{eqnarray}
which in turn leads to
\begin{eqnarray}
\left(R^{k+1}\right)^2 &=& \left(\left\|\vp^{k+1}-\vp^k\right\|^2 + \left\|\nabla\left(\vp^{k+1}-\vp^k\right)\right\|^2 + \left\|\ub^{k+1}-\ub^{k}\right\|^2\right)^2 \non\\
&\leq& C\left(\left\|\vp^{k+1}-\vp^k\right\|^4 + \left\|\nabla\left(\vp^{k+1}-\vp^k\right)\right\|^4 + \left\|\ub^{k+1}-\ub^{k}\right\|^4\right) \non\\
&\leq& C\tau^3 \int_{t_k}^{t_{k+1}}\left( \left\|\partial_t\vp(\cdot,t)\right\|^4 + \left\|\nabla\partial_t\vp(\cdot,t)\right\|^4 + \left\|\partial_t\ub(\cdot,t)\right\|^4\right) dt.
\end{eqnarray}
Therefore, for $d=2,3$, by using Cauchy-Schwarz inequality and Lemma \ref{numerical_estimate}, one gets
\begin{eqnarray}
&& \sum_{k=0}^{l}\left(1+\left\|\vp_h^k\right\|_{\infty}^2\right)R^{k+1} \leq \left(\sum_{k=0}^l C\left(1+\left\|\vp_h^k\right\|_{\infty}^{\frac{2(6-d)}{d}}\right)^{2} \right)^{1/2} \left(\sum_{k=0}^l\left(R^{k+1}\right)^2\right)^{1/2} \non\\
&\leq& \left(\sum_{k=0}^l C\left(1+\left\|\vp_h^k\right\|_{\infty}^{\frac{4(6-d)}{d}}\right) \right)^{1/2} \left(C\tau^3 \int_{0}^{t_{l+1}}\left( \left\|\partial_t\vp(\cdot,t)\right\|^4 + \left\|\nabla\partial_t\vp(\cdot,t)\right\|^4 + \left\|\partial_t\ub(\cdot,t)\right\|^4\right) dt\right)^{1/2} \non\\
&\leq& C\tau\left(\tau\sum_{k=0}^l\left(1+\left\|\vp_h^k\right\|_{\infty}^{\frac{4(6-d)}{d}}\right) \right)^{1/2} \left(\int_{0}^{t_{l+1}}\left( \left\|\partial_t\vp(\cdot,t)\right\|^4 + \left\|\nabla\partial_t\vp(\cdot,t)\right\|^4 + \left\|\partial_t\ub(\cdot,t)\right\|^4\right) dt\right)^{1/2} \non\\
&\leq& C\tau\sqrt{T+1} \left(\int_{0}^{t_{l+1}}\left( \left\|\partial_t\vp(\cdot,t)\right\|^4 + \left\|\nabla\partial_t\vp(\cdot,t)\right\|^4 + \left\|\partial_t\ub(\cdot,t)\right\|^4\right) dt\right)^{1/2} \leq C\tau\sqrt{T+1}.
\end{eqnarray}
Similarly, we have 
\begin{eqnarray}
\sum_{k=0}^{l}\left(1+\left\|\vp_h^{k+1}\right\|_{\infty}^4\right)\left\|\nabla\rho^{\vp,k+1}\right\|^2 &\leq& \tau^{-1/2}\left(\tau\sum_{k=0}^{l}C\left(1+\left\|\vp_h^{k+1}\right\|_{\infty}^{\frac{8(6-d)}{d}}\right)\right)^{1/2}
\left(\sum_{k=0}^{l}\left\|\nabla\rho^{\vp,k+1}\right\|^4\right)^{1/2} \non\\
&\leq& C\tau^{-1/2}\sqrt{T+1}\left(\sum_{k=0}^{l}\left\|\nabla\rho^{\vp,k+1}\right\|^4\right)^{1/2},
\end{eqnarray}
\begin{eqnarray}
\sum_{k=0}^{l}\left\|\vp_h^k\right\|_{\infty}^2 \left(\left\|\rho^{\ub,k}\right\|^2\right) &\leq& \tau^{-1/2}\left(\tau\sum_{k=0}^{l}\left(1 + \left\|\vp_h^k\right\|_{\infty}^{\frac{4(6-d)}{d}}\right)\right)^{1/2} \left(\sum_{k=0}^{l}2\left\|\rho^{\ub,k}\right\|^4\right)^{1/2} \non\\
&\leq& C\tau^{-1/2}\sqrt{T+1} \left(\sum_{k=0}^{l}\left\|\rho^{\ub,k}\right\|^4\right)^{1/2}.
\end{eqnarray}

Henceforth, it follows that 
\begin{eqnarray}
\sum_{k=0}^{l} \mathcal{R}^{k+1}
&=& \sum_{k=0}^l \bigg[
\tau^2 + \left\|R^{\vp,k+1}\right\|^2 + \left\|R_c^{\ub,k+1}\right\|^2 + \left\|R_m^{\ub,k+1}\right\|^2 + \left(1+\left\|\vp_h^k\right\|_{\infty}^2\right)R^{k+1} \non\\
&&\ \ \ \ \ \  + \left(1+\left\|\vp_h^{k+1}\right\|_{\infty}^4\right)\left\|\nabla\rho^{\vp,k+1}\right\|^2 + \left\|\nabla\rho^{\vp,k}\right\|^2 + \left\|\nabla\rho^{\mu,k+1}\right\|_6^2 + \left\|\vp_h^k\right\|_{\infty}^2 \left\|\rho^{\ub,k}\right\|^2\bigg] \non\\
&\leq& \left(T + C\sqrt{T+1}\right)\tau + \frac{2\tau}{3}\int_0^{t_{l+1}}\left(\left\|\partial_{tt}\vp(\cdot,t)\right\|^2 + \left\|\partial_{tt}\ub(\cdot,t)\right\|^2\right)dt \non\\
&&\ +\ \frac{2}{\tau}\int_{0}^{t_{l+1}}\left(\left\|\partial_t\rho^{\vp}(\cdot,t)\right\|^2 + \left\|\partial_t\rho^{\ub}(\cdot,t)\right\|^2\right)dt + \sum_{k=0}^l\left((1+\tau^2)\left\|\nabla\rho^{\mu,k+1}\right\|_6^2 + \left\|\nabla\rho^{\vp,k}\right\|^2\right) \non\\
&&\ +\ C\tau^{-1/2}\sqrt{T+1}\left(\sum_{k=0}^{l}\left\|\nabla\rho^{\vp,k+1}\right\|^4\right)^{1/2} + C\tau^{-1/2}\sqrt{T+1} \left(\sum_{k=0}^{l}\left\|\rho^{\ub,k}\right\|^4\right)^{1/2}. \non\\
&\leq& \left(T + C\sqrt{T+1} +\frac{2}{3}\right)\tau + \frac{2}{\tau}\int_{0}^{t_{l+1}}\left(\left\|\partial_t\rho^{\vp}(\cdot,t)\right\|^2 + \left\|\partial_t\rho^{\ub}(\cdot,t)\right\|^2\right)dt \non\\
&&\ +\ C\tau^{-1/2}\sqrt{T+1}\left(\sum_{k=0}^{l}\left\|\nabla\rho^{\vp,k+1}\right\|^4\right)^{1/2} + \sum_{k=0}^l\left((1+\tau^2)\left\|\nabla\rho^{\mu,k+1}\right\|_6^2 + \left\|\nabla\rho^{\vp,k}\right\|^2\right) \non\\
&&\ +\ C\tau^{-1/2}\sqrt{T+1} \left(\sum_{k=0}^{l}\left\|\rho^{\ub,k}\right\|^4\right)^{1/2}.
\end{eqnarray}
This completes the proof. \hfill 
\end{proof}

Now we are ready to prove  the main convergence theorem.
\begin{mytheo}\label{main_theo}
Suppose $(\varphi, \mu, \ub_c, \ub_m, P_c, P_m)$ is a weak solution to \eqref{ww_CH1}--\eqref{ww_D} with the additional regularities described in Assumption \ref{higher_regularities}. Recall the definition of  error functions $\sigma$s in Eqs. \eqref{ef_1}--\eqref{ef_4} and the $\rho^{\varphi}, \rho^{\mathbf{u}}, \rho^{\mu}$ in Eqs. \eqref{rho_1}--\eqref{rho_4}. Then, provided that $0<\tau<\tau_1$ for some sufficiently small $\tau_1>0$,
\begin{eqnarray}
&& \max_{0\leq k\leq K-1}\left(\left\|\nabla\sigma^{\varphi,k+1}\right\|^2 + \left\|\sigma_c^{\ub,k+1}\right\|^2 + \left\|\sigma_m^{\ub,k+1}\right\|^2\right) + \tau\sum_{k=0}^{K-1}\left\|\nabla\sigma^{\mu,k+1}\right\|^2 \nonumber\\
&& + \sum_{k=0}^{K-1}\left(\left\|\nabla(\sigma^{\varphi,k+1}-\sigma^{\varphi,k})\right\|^2  + \left\|\sigma_c^{\ub,k+1}-\sigma_c^{\ub,k}\right\|^2 + \left\|\sigma_m^{\ub,k+1}-\sigma_m^{\ub,k}\right\|^2 \right)\nonumber\\
&& + \tau\sum_{k=0}^{K-1}\left[\left\|\mathbb{D}(\sigma_c^{\ub,k+1})\right\|^2 + \sum_{i=1}^{d-1}\left\|\sigma_c^{\ub,k+1}\cdot\btau_i\right\|_{cm}^2 + \left\|\sigma_m^{\ub,k+1}\right\|^2\right] + \tau^2\sum_{k=0}^{K-1}\left\|\varphi_{h}^{k}\nabla\sigma^{\mu,k+1}\right\|^2 \nonumber\\
&\leq& C(T)\Bigg[\tau^2 + \int_{0}^{T}\left(\left\|\partial_t\rho^{\vp}(\cdot,t)\right\|^2 + \left\|\partial_t\rho^{\ub}(\cdot,t)\right\|^2\right)dt + \tau^{1/2}\bigg(\sum_{k=0}^{K}\left\|\nabla\rho^{\vp,k+1}\right\|^4\bigg)^{1/2} \non\\
&&\ \ \ \ \ \ \ \ \ \ +\ \tau\sum_{k=0}^K\left(\left\|\nabla\rho^{\vp,k}\right\|^2 + (1+\tau^2)\left\|\nabla\rho^{\mu,k+1}\right\|_6^2 \right) + \tau^{1/2}\bigg(\sum_{k=0}^{K}\left\|\rho^{\ub,k}\right\|^4\bigg)^{1/2} \Bigg] \label{final-err}
\end{eqnarray}
holds for some constant $C(T)>0$ independent of $\tau$ and $h$.
\end{mytheo}
\begin{proof}
Applying $\tau\sum_{k=0}^{l}$ to \eqref{combine}, and observing that $\sigma^{\varphi,k}\equiv 0$ and $\sigma_j^{\ub,k}\equiv 0$ for $k=0,\ j\in\{c,m\}$, it follows that 
\begin{eqnarray}
&& \frac{\gamma\epsilon}{2}\left\|\nabla\sigma^{\varphi,l+1}\right\|^2 + \frac{\rho_0}{2}\left\|\sigma_c^{\ub,l+1}\right\|^2 + \frac{\rho_0}{2\chi}\left\|\sigma_m^{\ub,l+1}\right\|^2 + \tau\sum_{k=0}^{l}\Big(\frac{M_0}{3}\|\nabla\sigma^{\mu,k+1}\|^2\Big) \nonumber\\
&& + \sum_{k=0}^{l}\left(\frac{\gamma\epsilon}{2}\left\|\nabla(\sigma^{\varphi,k+1}-\sigma^{\varphi,k})\right\|^2  + \frac{\rho_0}{2}\left\|\sigma_c^{\ub,k+1}-\sigma_c^{\ub,k}\right\|^2 + \frac{\rho_0}{2\chi}\left\|\sigma_m^{\ub,k+1}-\sigma_m^{\ub,k}\right\|^2 \right)\nonumber\\
&& + \tau\sum_{k=0}^{l}\left[\nu_0\left\|\mathbb{D}(\sigma_c^{\ub,k+1})\right\|^2 + \alpha_{BJSJ}\frac{\nu_0}{2\sqrt{d\lambda}}\sum_{i=1}^{d-1}\left\|\sigma_c^{\ub,k+1}\cdot\btau_i\right\|_{cm}^2 + \frac{\nu_0}{2\lambda}\left\|\sigma_m^{\ub,k+1}\right\|^2\right] \non\\
&& + \frac{\tau^2}{\rho_0}\sum_{k=0}^{l}\left[ \left\|\varphi_{c,h}^{k}\nabla\sigma_c^{\mu,k+1}\right\|^2 + \chi\left\|\varphi_{m,h}^{k}\nabla\sigma_m^{\mu,k+1}\right\|^2 \right] \nonumber\\
&\leq& \ C\tau\sum_{k=0}^l\mathcal{R}^{k+1} + C\tau\sum_{k=0}^l\left(1+\left\|\vp_h^{k+1}\right\|_{\infty}^{4}\right)\left\|\nabla\sigma^{\vp,k+1}\right\|^2 + \tau\sum_{k=0}^l\left(1+C\left\|\varphi_h^k\right\|_{\infty}^2\right) \left\|\sigma^{\ub,k+1}\right\|^2 \non\\
&& \ + C\tau\sum_{k=1}^l\left\|\nabla\sigma^{\vp,k}\right\|^2 + C\tau\sum_{k=1}^l\left\|\vp_h^k\right\|_{\infty}^2 \left\|\sigma^{\ub,k}\right\|^2 \non\\
&\leq& \ C\tau\sum_{k=0}^l\mathcal{R}^{k+1} + C\tau\left(1+\left\|\vp_h^{l+1}\right\|_{\infty}^{\frac{4(6-d)}{d}}\right)\left\|\nabla\sigma^{\vp,l+1}\right\|^2 + \tau\left(1+C\left\|\varphi_h^{l}\right\|_{\infty}^{\frac{2(6-d)}{d}}\right) \left\|\sigma^{\ub,l+1}\right\|^2 \non\\
&& \ +\ C\tau\sum_{k=1}^l\left(1+\left\|\vp_h^{k}\right\|_{\infty}^{\frac{4(6-d)}{d}}\right)\left\|\nabla\sigma^{\vp,k}\right\|^2 + C\tau\sum_{k=1}^l\left(1 + 2\left\|\vp_h^k\right\|_{\infty}^{\frac{2(6-d)}{d}}\right) \left\|\sigma^{\ub,k}\right\|^2.
\end{eqnarray}
Moving all the terms indexed  $(l+1)$ to the left hand side, one has
\begin{eqnarray}\label{sum_latest_step_set}
&& \Bigg(\frac{\gamma\epsilon}{2}-C\tau\Big(1+\big\|\varphi_h^{l+1}\big\|_{\infty}^{\frac{4(6-d)}{d}}\Big)\Bigg) \Big\|\nabla\sigma^{\varphi,l+1}\Big\|^2 + \Bigg(\frac{\rho_0}{2}-\tau\Big(1+C\left\|\varphi_h^l\right\|_{\infty}^{\frac{2(6-d)}{d}}\Big)\Bigg) \Big\|\sigma_c^{\ub,l+1}\Big\|^2 \nonumber\\
&& + \Bigg(\frac{\rho_0}{2\chi}-\tau\Big(1+C\left\|\varphi_h^l\right\|_{\infty}^{\frac{2(6-d)}{d}}\Big)\Bigg) \Big\|\sigma_m^{\ub,l+1}\Big\|^2 + \tau\sum_{k=0}^{l}\Big(\frac{M_0}{3}\left\|\nabla\sigma^{\mu,k+1}\right\|^2\Big) \nonumber\\
&& + \sum_{k=0}^{l}\left(\frac{\gamma\epsilon}{2}\Big\|\nabla(\sigma^{\varphi,k+1}-\sigma^{\varphi,k})\Big\|^2  + \frac{\rho_0}{2}\Big\|\sigma_c^{\ub,k+1}-\sigma_c^{\ub,k}\Big\|^2 + \frac{\rho_0}{2\chi}\Big\|\sigma_m^{\ub,k+1}-\sigma_m^{\ub,k}\Big\|^2 \right)\nonumber\\
&& + \tau\sum_{k=0}^{l}\left[\nu_0\left\|\mathbb{D}(\sigma_c^{\ub,k+1})\right\|^2 + \alpha_{BJSJ}\frac{\nu_0}{2\sqrt{d\lambda}}\sum_{i=1}^{d-1}\left\|\sigma_c^{\ub,k+1}\cdot\btau_i\right\|_{cm}^2 + \frac{\nu_0}{2\lambda}\left\|\sigma_m^{\ub,k+1}\right\|^2\right] \non\\
&& + \frac{\tau^2}{\rho_0}\sum_{k=0}^{l}\left[ \left\|\varphi_{c,h}^{k}\nabla\sigma_c^{\mu,k+1}\right\|^2 + \chi\left\|\varphi_{m,h}^{k}\nabla\sigma_m^{\mu,k+1}\right\|^2 \right] \nonumber\\
&\leq& C\tau\sum_{k=0}^l \mathcal{R}^{k+1} + C\tau\sum_{k=1}^l\left(1+\left\|\vp_h^{k}\right\|_{\infty}^{\frac{4(6-d)}{d}}\right)\left\|\nabla\sigma^{\vp,k}\right\|^2 + C\tau\sum_{k=1}^l\left(1 + 2\left\|\vp_h^k\right\|_{\infty}^{\frac{2(6-d)}{d}}\right) \left\|\sigma^{\ub,k}\right\|^2.
\end{eqnarray}
By Lemma \ref{numerical_estimate} we have, for all $0\leq l \leq K-1$, 
\begin{eqnarray}
\tau^{\frac{1}{2}}\left\|\varphi_h^{l+1}\right\|_{\infty}^{\frac{4(6-d)}{d}} =
\left(\tau\left\|\varphi_h^{l+1}\right\|_{\infty}^{\frac{8(6-d)}{d}}\right)^{\frac{1}{2}} \leq \left(\tau\sum_{k=0}^{K-1}\left\|\varphi_h^{k+1}\right\|_{\infty}^{\frac{8(6-d)}{d}}\right)^{\frac{1}{2}} \leq C\sqrt{T+1}.
\end{eqnarray}
Hence we can choose a sufficiently small $\tau_1$ such that for all $0<\tau<\tau_1$ and $0\leq l \leq K-1$
\begin{eqnarray}
&&C\tau\left(1+\left\|\varphi_h^{l+1}\right\|_{\infty}^{\frac{4(6-d)}{d}}\right) \leq C\tau + C\tau^{\frac{1}{2}}\Big(C\sqrt{T+1}\Big) \leq \frac{\gamma\epsilon}{4}, \\
&& \frac{\gamma\epsilon}{2}-C\tau\bigg(1+\left\|\varphi_h^{l+1}\right\|_{\infty}^{\frac{4(6-d)}{d}}\bigg) \geq \frac{\gamma\epsilon}{4}\ , \\
&& \frac{\rho_0}{2}-C\tau\bigg(1+\left\|\varphi_{h}^l\right\|_{\infty}^{\frac{2(6-d)}{d}}\bigg) \geq \frac{\rho_0}{4}\ , \\
&& \frac{\rho_0}{2\chi}-C\tau\bigg(1+\left\|\varphi_{h}^l\right\|_{\infty}^{\frac{2(6-d)}{d}}\bigg) \geq \frac{\rho_0}{4\chi}\ .
\end{eqnarray}
It follows from \eqref{sum_latest_step_set} that
\begin{eqnarray}
&& \frac{\gamma\epsilon}{4} \Big\|\nabla\sigma^{\varphi,l+1}\Big\|^2 + \frac{\rho_0}{4}\Big\|\sigma_c^{\ub,l+1}\Big\|^2 + \frac{\rho_0}{4\chi}\Big\|\sigma_m^{\ub,l+1}\Big\|^2 + \tau\sum_{k=0}^{l}\Big(\frac{M_0}{3}\Big\|\nabla\sigma^{\mu,k+1}\Big\|^2\Big) \nonumber\\
&& + \sum_{k=0}^{l}\left(\frac{\gamma\epsilon}{2}\Big\|\nabla(\sigma^{\varphi,k+1}-\sigma^{\varphi,k})\Big\|^2  + \frac{\rho_0}{2}\Big\|\sigma_c^{\ub,k+1}-\sigma_c^{\ub,k}\Big\|^2 + \frac{\rho_0}{2\chi}\Big\|\sigma_m^{\ub,k+1}-\sigma_m^{\ub,k}\Big\|^2 \right)\nonumber\\
&& + \tau\sum_{k=0}^{l}\left[\nu_0\left\|\mathbb{D}(\sigma_c^{\ub,k+1})\right\|^2 + \alpha_{BJSJ}\frac{\nu_0}{2\sqrt{d\lambda}}\sum_{i=1}^{d-1}\left\|\sigma_c^{\ub,k+1}\cdot\btau_i\right\|_{cm}^2 + \frac{\nu_0}{2\lambda}\left\|\sigma_m^{\ub,k+1}\right\|^2\right] \non\\
&& + \frac{\tau^2}{\rho_0}\sum_{k=0}^{l}\left[ \left\|\varphi_{c,h}^{k}\nabla\sigma_c^{\mu,k+1}\right\|^2 + \chi\left\|\varphi_{m,h}^{k}\nabla\sigma_m^{\mu,k+1}\right\|^2 \right] \nonumber\\
&\leq& C\tau\sum_{k=0}^l \mathcal{R}^{k+1} + C\tau\sum_{k=1}^l\left(1+\left\|\vp_h^{k}\right\|_{\infty}^{\frac{4(6-d)}{d}}\right)\left\|\nabla\sigma^{\vp,k}\right\|^2 + C\tau\sum_{k=1}^l\left(1 + 2\left\|\vp_h^k\right\|_{\infty}^{\frac{2(6-d)}{d}}\right) \left\|\sigma^{\ub,k}\right\|^2. \nonumber 
\end{eqnarray}
 Noticing that  $\tau\sum_{k=0}^{K}\Big\|\varphi_h^{k}\Big\|_{\infty}^{\frac{p(6-d)}{d}} \leq C(T+1)$ for $p=2,4$  and in light of  Lemma \ref{sum-R}, we arrive at  the error estimate \eqref{final-err}  by setting $l=K-1$ and applying discrete Gronwall's inequality. This completes the proof.  
\end{proof}

\begin{corollary}\label{main_cor}
	Suppose $(\varphi, \mu, \ub_c, \ub_m, P_c, P_m)$ is a weak solution to \eqref{ww_CH1}--\eqref{ww_D} satisfying the  regularities  Assumption \ref{higher_regularities}. Then there exists  $\tau_1>0$ such that for all $\tau<\tau_1$  the following optimal convergence rates hold
	\begin{eqnarray*}
		&& \max_{0\leq k\leq K-1}\left(\left\|\nabla e^{\varphi,k+1}\right\|^2 + \left\|e_c^{\ub,k+1}\right\|^2 + \left\|e_m^{\ub,k+1}\right\|^2\right) + \tau\sum_{k=0}^{K-1}\left\|\nabla e^{\mu,k+1}\right\|^2  + \tau\sum_{k=0}^{K-1}\left\|\mathbb{D}(e_c^{\ub,k+1})\right\|^2 \\
		&&\leq C(T)\big(\tau^2 +h^{2q}\big),
	\end{eqnarray*}
	where $q\geq 1$ is the spatial approximation order.
\end{corollary}
For numerical evidence of the convergence results, we refer to \cite{CHW2017}.

\begin{remark}
	In the the discrete energy dissipation analysis established in Chen et al. (2017), for the numerical scheme, a cancellation of a nonlinear error term associated with the convection part has played a very important role. Meanwhile, in the optimal rate error estimate presented in this section, such a cancellation technique is not needed in the convergence proof, due to the subtle fact that, a growth constant for the velocity error term, namely $(1 + C \|\varphi_h^{k} \|_{\infty}^2 )$ appearing in (4.49), would not lead to a theoretical difficulty in the derivation of discrete Gronwall inequality. This fact is associated with Nativer-Stokes nature for the fluid velocity, in which the higher order kinematic diffusion and the temporal derivative of the velocity variable have greatly facilitated the analysis at both the analytic and numerical levels. In comparison, for the Cahn-Hillird-Hele-Shaw system, in which the fluid velocity is statically determined by the phase field variables, such a cancellation technique is necessary to pass through the optimal rate convergence analysis, because of lack of regularity for the velocity field; see the related works Chen et al. (2016); Diegel et al. (2017); Liu et al. (2017),  etc.
\end{remark}

\section{Concluding remarks}\label{sec:conclusion}

In this article we provide an optimal rate convergence analysis and error estimate of a fully discrete finite element numerical scheme for the Cahn-Hilliard-Stokes-Darcy system that models two-phase flows. An operator splitting is applied in the numerical scheme, so that a coupling between the Cahn-Hilliard and the fluid solvers is avoided. The unique solvability and the energy stability have already been proved in the existing literature. The optimal rate error estimate is established in the energy norm, $\ell^\infty (0, T; H^1) \cap \ell^2 (0, T; H^2)$ norm for the phase variables, and $\ell^\infty (0, T; H^1) \cap \ell^2 (0, T; H^2)$ norm for the velocity variable.  A discrete $\ell^2 (0;T; H^3)$ bound of the numerical solution for the phase variables also plays an important role, which is accomplished via a discrete version of Gagliardo-Nirenberg inequality in the finite element space. 
	
\section*{Acknowledgements}
W. Chen is supported by the National Key R\&D Program of China (2019YFA0709502) and the National Science Foundation of China (12071090).  D. Han acknowledges support from NSF-DMS-1912715. C. Wang is supported by NSF DMS-2012669.  C.Wang also thanks the Key Laboratory of Mathematics for Nonlinear Sciences, Fudan University for support during his visit. X. Wang thanks support from NSFC11871159, Guangdong Provincial Key Laboratory of Computational Science and Material Design via 2019B030301001.

\bibliographystyle{plain}
\bibliography{reference_err,multiphase-2020}

\end{document}